\documentclass[a4paper,11pt]{article}

\usepackage{amsmath,amsfonts,latexsym,theorem}
\usepackage{pdfpages}
\usepackage{subfigure}
\usepackage{caption}
\nonstopmode
\usepackage{graphicx}
\usepackage{epsfig}
\usepackage{epstopdf}

\numberwithin{equation}{section}

\newtheorem{Theorem}{Theorem}[section]
\newtheorem{Definition}{Definition}[section]
\newtheorem{Proposition}{Proposition}[section]
\newtheorem{Lemma}{Lemma}[section]

\newenvironment{Proofc}[1]{\smallskip\par\noindent\textsc{#1}\quad}%
  {\hfill$\Box$\bigskip\par}
\newenvironment{Proof}{\begin{Proofc}{Proof}}{\end{Proofc}}

\newtheorem{Remark}{Remark}[section]

\def\a{\alpha}
\def\b{\beta}
\def\d{\delta}

\def\g{\gamma}
\def\G{\Gamma}

\def\r{\rho}
\def\s{\sigma}
\def\o{\omega}

\def\gS{\Sigma}

\def\pd{\partial}

\newcommand{\cR}{{\cal R}}

\newcommand{\R}{{\mathbb R}}
\newcommand{\N}{{\mathbb N}}

\def\pd{\partial}

\begin{document}
\title{A  flame propagation model  on a  network with application to  a blocking problem}

\author{Fabio Camilli\footnotemark[1] \and Elisabetta Carlini \footnotemark[1] \and Claudio Marchi \footnotemark[2] }
\date{version: \today}
\maketitle

\footnotetext[1]{Dip. di Scienze di Base e Applicate per l'Ingegneria,  ``Sapienza'' Universit{\`a}  di Roma, via Scarpa 16,
 00161 Roma, Italy, ({\tt e-mail: camilli@sbai.uniroma1.it})}
 \footnotetext[2]{Dip. di Matematica,  ``Sapienza'' Universit{\`a}  di Roma, p.le A.Moro 2,
 00161 Roma, Italy, ({\tt e-mail: carlini@mat.uniroma1.it})}
\footnotetext[3]{Dip. di Ingegneria dell'Informazione, Universit{\`a} di Padova, via Gradenigo 6/B, 35123 Padova, Italy ({\tt claudio.marchi@unipd.it}).}

\begin{abstract}
We consider the Cauchy problem
\[
    \left\{\begin{array}{ll}
       \pd_t u+H(x,Du)=0  \quad &(x,t)\in\G\times (0,T),\\
       u(x,0)=u_0(x)  &x\in\G
     \end{array}
     \right.
\]
where $\G$ is a network and $H$ is a convex and  positive homogeneous Hamiltonian which may change from edge to edge. In the former part of the paper, we prove that the Hopf-Lax type formula gives the (unique) viscosity solution of the problem. In the latter part of the paper we study  a flame propagation model in a network and an optimal strategy to block a fire breaking up in some part of a pipeline; some numerical simulations are provided.
\end{abstract}
 \begin{description}
\item [{\bf MSC 2000}:]  49L25, 35F21, 65M06
\item [{\bf Keywords}:]    Evolutive Hamilton-Jacobi equation, viscosity solution, network, Hopf-Lax formula, approximation.
\end{description}
\section{Introduction}
We study the Cauchy problem
\begin{equation}\label{int-Cauchy}
    \left\{\begin{array}{ll}
       \pd_t u+H(x,Du)=0  \quad &(x,t)\in\G\times (0,T)\\
       u(x,0)=u_0(x)  &x\in\G
     \end{array}
     \right.
\end{equation}
where $\G$ is a network and the operator $H : \G\times \R \to\R$ may change from edge to edge and inside each edge it is  continuous, convex, nonnegative and positive homogeneous in the last variable.\par
When the state variable varies in an Euclidean space~$\R^n$, the problem \eqref{int-Cauchy} arises in flame propagation models and evolution of curves whose speed of propagation only depends on the normal direction. Existence, uniqueness and evolution of level sets of the solution of \eqref{int-Cauchy} have been extensively studied in the framework of viscosity solution theory (see \cite{Ba,Bass,cgg,so}). In this case, the unique viscosity solution
of \eqref{int-Cauchy} is given by the Hopf-Lax formula
\begin{equation}\label{int-hopf}
u(x, t) =\min\{u_0(y):\, S(y, x)\le t\},\qquad (x,t)\in\R^n\times (0,\infty),
\end{equation}
where $S$ is a distance function characterized by solving the associated stationary equation
  \[ H(x,Dw) = 1,\qquad x\in\R^n.\]
In the recent time, there is an increasing interest in the study of nonlinear differential equations on networks since they
describe  various phenomena as traffic flow, blood circulation, data transmission, electric networks, etc (see \cite{gp,m}). Concerning Hamilton-Jacobi equations  on networks, we mention the recent papers \cite{acct,IM,ls,sc} where different notions of viscosity solution  have been introduced; we refer to \cite{CMcompar} for a comparison among some of  them.

This paper is divided into two parts; in the former one, following the approach in \cite{ls},  we 
prove that the Hopf-Lax formula \eqref{int-hopf} can be extended to this framework.
The main issue of the investigation is to tackle transition vertices (namely, points of the network where several edges meet each other). Actually,   a suitable definition of viscosity solution at transition vertices (together with the standard one at points inside edges) will ensure the well posedness of the problem. Let us recall that this feature also happens for stationary first order equations (see \cite{acct,IM,ls,sc}) whereas, for second order equations, some transition conditions (the so-called Kirchhoff condition) need to be imposed (see \cite{cms,m} and references therein).\par
In the second part of the paper we illustrate our results with a concrete application: the {\it blocking} problem. Suppose that a fire breaks up in some part of an  oil pipeline.
 A central controller can   stop the propagation of the fire  by closing the junctions of the pipes, represented by the vertices of the network. The controller spends  some time   to reach the junctions which   become unavailable when they  are reached by the fire front. Therefore only a subset of the vertices can be closed on time  to stop the fire. The aim is to find a strategy which maximizes the part of the network preserved by the fire.
 We give a characterization of  the optimal strategy and we study the corresponding flame propagation in the the network. Moreover we describe  a numerical scheme for the solution of the problem  and we present  some numerical examples.

This paper is organized as follows: in the rest of the Introduction we  set our notations. Section~\ref{sect:2} is devoted to the theoretical problem; 
Section~\ref{sect:3} is devoted to the application to the blocking problem also providing some numerical simulations.

\medskip
\noindent\textbf{Notations:}
A  network  $\G$ is a connected subset of $\R^n$ formed by a finite collection of points $V:=\{x_i\}_{i\in I}$ and edges $E:=\{e_j\}_{j\in J}$. The vertices of $V$ are connected by  the continuous, non self-intersecting arcs of $E$. Each arc $e_j$ is  parametrized by a smooth function $\pi_j:[0,l_j]\to\R^n,\, l_j>0$ and we set
\[
e_j=\pi_j((0,l_j))\qquad\hbox{and}\qquad \overline e_j=\pi_j([0,l_j])\,.
\]
For $i\in I$ we denote by $Inc_i:=\{j\in J\mid\,e_j \,\text{is incident to}\,x_i\}$ the set of arcs incident a same vertex $x_i$. We  fix  a set   $I_B\subset I$  and we denote  by  $\partial \G:=\{x_i\in V\mid \,i\in I_B\}$, the set of boundary vertices of $\G$.
We denote by $ d:\G\times\G\to \R^+$ the path distance on $\G$, i.e.
 \begin{equation}\label{geodistC}
   d(x,y):=\inf\left\{\ell(\g):\,\text{$\g\subset \G$ is a path joining $x$ to $y$ }\right\}, \qquad x,y\in \G
 \end{equation}
where $\ell(\g)$ is the length of $\g$. We assume that the network is connected, hence $d(x,y)$ is finite for any $x,y\in \G$. \\
We shall always identify   $x\in  \overline e_j$ with   $y=\pi_j^{-1}(x)\in [0,l_j]$. For any function $u:  \G\to\R$ and each $j\in J$ we denote by $u_j:[0,l_j]\to \R$ the restriction of $u$ to $\overline e_j$, i.e.  $u_j(y)=u(\pi_j(y))$ for $y\in [0,l_j]$. The derivative are always considered  with respect to the parametrization of the arc, i.e. if $x\in e_j$,  $y=\pi^{-1}_j(x)$ then
$Du(x):= \frac{d  u_j}{d y }  (y)$. At $x=x_i\in V$, we denote $D_ju(x_i)$ the internal derivative relative to the arc $e_j$, $ j\in Inc_i$, i.e.
\begin{equation*}
D_j u (x_i)=\left\{
               \begin{array}{ll}
                 \lim\limits_{h\to 0^+} \frac{u_j(h)-u_j(0)}{h}, & \hbox{if $x_i=\pi_j(0)$,} \\[4pt]
                 \lim\limits_{h\to 0^+}\frac{u_j( l_j-h)-u_j( l_j)}{h}, & \hbox{if $x_i=\pi_j( l_j)$.}
               \end{array}
             \right.
\end{equation*}
The space $C(\G\times (0,T))$ of continuous functions on $\G\times (0,T)$ is the space of $u:\G\times (0,T)\to\R$ such that
$u_j \in C([0, l_j]\times (0,T))$ and $u_j(x_i,t)=u_k(x_i,t)$ for all $j,k\in Inc_i$, $t\in(0,T)$ and all $i\in I$.

\section{Evolutive Hamilton-Jacobi equations on networks}\label{sect:2}
In this section we assume for simplicity  that $\pd \G=\emptyset$ (otherwise it is possible to introduce appropriate boundary condition
on $\pd \G$) and we consider the Hamilton-Jacobi equation
\begin{equation}\label{HJe}
    \pd_t u+H(x,Du)=0\qquad  (x,t)\in \G\times (0,T)
\end{equation}
with   the initial condition
\begin{equation}\label{IC}
  u(x,0)=u_0(x_0),\qquad x\in \G.
\end{equation}
The   Hamiltonian $H$ is given by  a  family $\{H_j\}_{j\in J}$, where   $H_j:\overline e_j\times \R\to\R$,  satisfies the following assumptions
\begin{align}
     &H_j\in C(\overline e_j\times  (0,T) );\label{H0}\\
    &|H_j(x,p)-H_j(y,p)|\le Cd(x,y)(1+|p|)\quad\text{for any $x$, $y\in \overline e_j$, $p\in \R$};\label{H1}\\
    &H_j(x,\cdot) \quad \text{is convex and  positive homogeneous  in $p$ for any $x\in \overline e_j$};\label{H2}\\
    &\inf\{H_j(x,p):\,|p|=1,\,x\in \overline e_j\}>0.  \label{H3}
    \end{align}
 \begin{Remark}
By  \eqref{H2} the equation is geometric and it is connected with front propagation (see \cite{Ba,Bass,so}). Assumption \eqref{H2} also implies the coercivity of the Hamiltonian, i.e. $\lim_{|p|\to\infty} H(x,p)=+\infty$ for any $x\in\G$.\\
A Hamiltonian satisfying the previous assumptions is given by
\[H_j(x,p)=\sup_{a\in A_j}\{-b_j(x,a)p\}\]
where $A_j$ is a compact metric space,  $b_j:\overline e_j\times A_j\to\R$ is a continuous function such that, for some $r>0$, there holds $(-r,r)\subset \overline{\text{co}}\{b_j(x,a):\,a\in A_j\}$.
In particular,  if $A_j=[-1,1]$ and $b_j(x,a)=c_j(x)a$ with $c$  bounded and strictly positive, then $H_j(x,p)=c_j(x)|p|$.
\end{Remark}
Let us now recall the definition of viscosity solution introduced in \cite{ls}. We consider the following class of test functions
\[
C^1(\G\times (0,T)):=\left\{\phi\in C(\G\times (0,T))\mid \phi_j\in C^1([0, l_j]\times (0,T))\quad \forall j\in J\right\}.
\]
\begin{Definition}\label{def_visco}\mbox{}
\begin{itemize}
\item[i)]  A function $u\in C(\G\times(0,T))$ is said a (viscosity) subsolution to~\eqref{HJe}
    if for every test function~$\phi\in C^1(\G\times(0,T))$ such that $u-\phi$ attains a local maximum   at $(x,t)\in  e_j \times(0,T)$, we have
  \begin{equation}\label{subsolv}
   \pd_t  \phi_j (x,t)+H_j(x, D\phi_j(x,t))   \leq 0.
  \end{equation}
\item [ii)] A function $u\in C( \G\times(0,T))$ is said a (viscosity) supersolution to~\eqref{HJe} if for every test function~$\phi\in C^1(\G\times(0,T))$ such that $u-\phi$ attains a local minimum   at $(x,t)\in\G\times(0,T)$, we have

\begin{equation}\label{supersolv}
\begin{split}
\pd_t  \phi (x,t)+H(x, D\phi(x,t))\ge 0 \qquad & \textrm{if }x\notin V\\[4pt]
\max_{j\in Inc_i}\{\pd_t  \phi_j (x,t)+H_j(x, D_j\phi(x,t))\}\ge 0  \qquad & \textrm{if }x=x_i\in V.
\end{split}
\end{equation}
\end{itemize}
A function $u\in C(\G\times(0,T) )$ is said  a (viscosity) solution of  \eqref{HJe}  if it is both a viscosity subsolution and a viscosity supersolution of \eqref{HJe}.
\end{Definition}
\begin{Remark}\label{rmk:1}
We observe that the continuity of~$\phi$ at $x_i\in V$ implies: $\phi_j(x_i,t)=\phi_k(x_i,t)$ for every $t\in[0,T]$ and $j,k\in Inc_i$. Hence, the transition condition in \eqref{supersolv} is equivalent to
\[
\pd_t  \phi (x_i,t)+\max_{j\in Inc_i}\{H_j(x_i, D_j\phi(x_i,t))\}\ge 0.
\]
\end{Remark}
 We first collect several properties of the viscosity solution of \eqref{HJe}.
The first result is  a comparison principle for the equation \eqref{HJe} established in~\cite{ls}.  
\begin{Proposition}\label{th:comparison}
Let $u,v\in C(\G\times [0,T])$  be a subsolution and, respectively, a supersolution of \eqref{HJe}
such that $u(x,0)\le v(x,0)$ for $x\in\G$. Then $u\le v$ in $\G\times [0,T]$.
\end{Proposition}
In the next Proposition~\ref{prp:reg} we state a regularity result for the solution of \eqref{HJe}.
\begin{Proposition}\label{prp:reg}
Let $u$ be a solution to \eqref{HJe}-\eqref{IC} where $u_0$ is Lipschitz continuous. Then, $u$ is Lipschitz continuous in $\G\times[0,T]$.
\end{Proposition}
The proof of the last proposition is a standard adaptation of the one in the Euclidean case (see \cite[Prop.2.1]{nr}) and we skip it. \par
We now  exploit the geometric character of the Hamiltonian (see \eqref{H2}) to give a representation formula for the solution of \eqref{HJe}. To this end, it is expedient to introduce some notations.
Given the Hamiltonian $H=\{H_j\}_{j\in J}$,  we define the support function of the sub-level set $\{p\in \R:\,H_j(x,p)\le 1\}$ as
\[ s_j(x,q):=\sup\{p\,q:\,H_j(x,p)\le1\} \]
and we set $ s:=\{ s_j\}_{j\in J}$.
The function  $s_j: \bar e_j\times \R\to\R$ is continuous, convex, positive homogeneous and non negative (see   \cite{s2}).
For example, if $H_j(x,p)=c_j(x)|p|$ then $s_j(x,q)=|q|/c_j(x)$.\par
A path $\xi:[0,t]\to\G$ is said {\it admissible} if  there are $t_0=0<t_1<\dots<t_{M+1}=t$ such that, for any $m=0,\dots,M$,   $\xi([t_m,t_{m+1}])\subset  \overline e_{j_m}$ for some $j_m\in J$ and $\pi_{j_m}^{-1}\circ \xi\in C^1(t_m,t_{m+1})$. We denote by $B^t_{y,x}$ the set of the admissible path such that $\xi(0)=y$, $\xi(t)=x$. We  introduce  a distance function related to the Hamiltonian $H$ on the network. For $x,y\in\G$ define
\begin{equation}\label{dist}
   S(y,x)=\inf\left\{\int_0^t  s(\xi(r),\dot\xi(r))dr:\,t>0,\,\xi\in B^t_{y,x}\right\},
\end{equation}
Note that the  distance defined by \eqref{dist} coincides with  the one  defined by \eqref{geodistC}   for  $H_j(x,p)=|p|$ for every $j\in J$; actually, in this case, we have $s_j(x,q)=|q|$  for every $j\in J$.
The next proposition summarizes some properties of $S$ (for the definition of viscosity solution on a network in the stationary case, we refer the reader to the paper \cite{ls}). 
\begin{Proposition}\label{prp:dist}
$S$ is a  Lipschitz continuous function on $\G\times \G$ and it is equivalent to the path distance $d$, i.e. there exists $C>0$ such that
\begin{equation}\label{equivdist}
   Cd(x,y)\le S(x,y)\le \frac{1}{C}\,d  (x,y),\qquad \text{for any $x,y\in\G$.}
\end{equation}
 Moreover, for any $K\subset\G$  closed,    $S(K,\cdot) $ is a subsolution in $\G$ and a supersolution in $\G\setminus K$ of the Hamilton-Jacobi equation
\begin{equation}\label{HJ}
    H(x,Du)=1.
\end{equation}
\end{Proposition}
\begin{Proof}
For the proof of the properties of $S$ and of bounds~\eqref{equivdist}   we refer to \cite[Prop.4.1]{sc} since their arguments easily adapt to our case.
We first prove that, for a given $x_0\in \G$,  $u(\cdot)=S(x_0,\cdot)$ is a viscosity subsolution in $\G$ and a viscosity
supersolution in $\G\setminus\{x_0\}$ of \eqref{HJ} in the sense of~\cite{ls}. This amount to prove that: $(i)$ $u$ is a viscosity subsolution  in $\G \setminus V$ and a viscosity supersolution in $(\G\setminus \{x_0\}) \setminus V$ in the standard viscosity solution sense (see \cite{bcd}), $(ii)$ for any $x=x_i\in V\setminus\{x_0\}$ and for any test function $\phi\in
C^1(\G):=\left\{\phi\in C(\G)\mid \phi_j\in C^1([0, l_j])\quad \forall j\in J\right\}$ such that $u-\phi$ has a local minimum
point at $x$, then
\begin{equation}\label{HJsuper}
    \max_{j\in Inc_i}\{ H_j(x, D_j\phi(x))\}\ge 1.
\end{equation}
For the proof that $S$ is a viscosity solution of \eqref{HJ} inside the edges  we refer to \cite[Thm.2.1]{s2}. We show the  the function $u$  satisfies \eqref{HJsuper} at a vertex $x=x_i\in V\setminus\{x_0\}$. Hence let $\phi$ be a test function such that $u-\phi$ has a local minimum at $x$ and let $\xi\in B^t_{x_0,x}$ be such that $u(x)=\int_0^t s( \xi, \dot\xi)ds$ (the coercivity of~$s$ and the $1$-dimensionality of the problem ensure the existence of such an optimal path; moreover, without any loss of generality, we assume that~$\xi$ is smooth). By definition of $\xi\in B^t_{x_0,x}$,  we can find  $\eta>0$ and  $j\in Inc_i$ such that $\xi(r)\in e_j$
 for $r\in (t-\eta,t)$. Hence
 \begin{align*}
    &u(\xi(r))=u_j(\pi_j^{-1}(\xi(r))),\quad \phi(\xi(r))=\phi_j(\pi_j^{-1}(\xi(r))) \qquad\text{for $r\in (t-\eta,t)$, }\\
    &u(x)=u_j(\pi_j^{-1}(\xi(t))),\quad \phi(x)=\phi_j(\pi_j^{-1}(\xi(t))).
 \end{align*}
 Since $u-\phi$ has a minimum point at $x$, we have
\begin{align*}
    \frac{\phi_j(\pi_j^{-1}(\xi(t)))-\phi_j(\pi_j^{-1}(\xi(r)))}{t-r}\ge \frac{u_j(\pi_j^{-1}(\xi(t)))-u_j(\pi_j^{-1}(\xi(r)))}{t-r}=\\
    \frac{1}{t-r}\int_r^t  s_j\left(\pi_j^{-1}(\xi(\r)),\frac{d}{d\r}(\pi_j^{-1}\circ \xi)(\r)\right)d\r.
\end{align*}
Passing to the limit for $r\to t$ we get
\[-D_j\phi(x)q\ge  s_j(\pi_j^{-1}(\xi(t)),-q)\]
where  $q= \frac{d}{d\r}(\pi_j^{-1}\circ \xi)(t)$. Hence, recalling the definition of $s_j$, it follows that
\[
  H_j(x,D_j\phi)\ge 1
 \]
and therefore \eqref{HJsuper}. \\
Having proved that $S(x_0,\cdot)$ is a viscosity subsolution in $\G$ and a viscosity
supersolution in $\G\setminus\{x_0\}$, then it is easy to prove that    $S(K,\cdot)$   is a subsolution in $\G$ and a supersolution in $\G\setminus K$ of \eqref{HJ} (see \cite[Prop.6.2]{sc}).
\end{Proof}
The following result gives a representation formula of Hopf-Lax type for the solution of \eqref{HJe}-\eqref{IC}. 
\begin{Theorem}\label{existenceC}
Let $u_0:\G\to\R$ be a continuous function. Then the solution of \eqref{HJe}-\eqref{IC}
 is given by
\begin{equation}\label{31}
  u(x,t)=\min\{u_0(y):\, S(y,x)\le t\}.
\end{equation}
\end{Theorem}
In order to prove this result, let us first establish some  preliminary lemmas.
\begin{Lemma}\label{caratstaz}
If $w$ is a subsolution (resp., supersolution) of \eqref{HJ} in $\G$, then $u(x,t)=w(x)-t$  is a subsolution (resp., supersolution) of \eqref{HJe} in $\G\times (0,T)$.
\end{Lemma}
\begin{Proof}
We only  show that $u$ is a supersolution at $(x_0,t_0)$ in the case  $x_0=x_i\in V$,  being the case $x_0\not\in V$ similar. Assume by contradiction that there exists  $\eta>0$ and a  test function $\phi$  at $(x_0,t_0)$ such that  $u-\phi$ has a local minimum at $(x_0,t_0)$
  and
\begin{equation}\label{eq1}
\max_{j\in Inc_i}\{ \pd_t \phi_j(x_0,t_0)+H_j(x_0, D_j\phi(x_0,t_0))\}\le -\eta<0.
\end{equation}
By
\begin{equation}\label{eq2}
 u(x_0,t_0)-\phi(x_0,t_0)\le u(x ,t)-\phi(x,t)
\end{equation}
for $x=x_0$,  taking into account the definition of $u$, we get
$\phi(x_0,t)-\phi(x_0,t_0)\le t_0-t$ $\forall t\in(0,T)$; by the arbitrariness of $t$, we infer
\begin{equation}\label{eq3}
\pd_t \phi_j(x_0,t_0)=-1,\qquad \forall j\in Inc_i.
\end{equation}
 Moreover by \eqref{eq2} for $t=t_0$ we get that $w(x)-\phi(x,t_0)$ has a local minimum at $x_0$, hence
 \begin{equation}\label{eq4}
 \max_{j\in Inc_i} \{H_j(x_0, D\phi(x_0,t_0))\}\ge 1.
 \end{equation}
By \eqref{eq3}, \eqref{eq4} and Remark~\ref{rmk:1}, we get a contradiction to \eqref{eq1}.
\end{Proof}


\begin{Lemma}\label{geometric}\hfill
\begin{itemize}
\item[$(i)$]
  A function $u\in C^0(\G  \times(0,T))$ is a subsolution of \eqref{HJe}  if and only if for any $\a\in\R$ and for any  admissible test function $\phi$  which has a local minimum on $\{u\ge \a\}\cap   (\G\times (0,T))$ at  $(x,t)\in e_j\times (0,T)$, then \eqref{subsolv} holds.
\item[$(ii)$]
A function $v\in C^0(\G\times (0,T))$ is a supersolution of \eqref{HJe}
if and only if for any $\a\in\R$ and for any  admissible test function $\phi$  which has a local maximum on $\{v\le \a\}\cap   (\G\times (0,T))$ at  $(x,t)\in \G\times (0,T)$, then \eqref{supersolv} holds.
\end{itemize}
\end{Lemma}
The proof of the previous lemma is similar to the one in \cite[Lemma3.1]{s1} which can be easily adapted to networks.

\begin{Proofc}{Proof of Theorem \ref{existenceC}}
We first prove that $u$ is continuous.
Given $(x_0,t_0)\in\G\times [0,T]$, let $(x_n,t_n)\in\G\times [0,T]$ be such that $\lim_{n\to\infty} (x_n,t_n)=(x_0,t_0)$ and set $\d_n=|t_n-t_0|+Cd (x_n,x_0)$ where $C$ is as in \eqref{equivdist}. We claim that
\begin{equation}\label{eq7}
\left\{y \in \G:\,S(y,x_0)\le t_0\right\}\subset\left\{y \in \G  :\,S(y,x_n)\le
t_n+\d_n\right\}
\end{equation}
In fact, if $S(y,x_0)\le t_0$, then, by \eqref{equivdist}, we get
\[S(y,x_n)\le S(y,x_0)+S(x_0,x_n)\le t_0+Cd(x_0,x_n)\le  t_n+\d_n\]
and therefore \eqref{eq7}.
Moreover
\begin{equation}\label{eq8}
    \left\{y \in \G:\,S(y,x_n)\le t_n+\d_n\right\}\subset \big\{y\in\G:\,d  \left(y,\{z:\,S(z,x_n)\le t_n\}\right)\le \d_n/C\big\}.
\end{equation}
In fact if $S(y,x_n)\le t_n+\d_n$, then
\begin{align*}
d \left(y,\{z:\,S(z,x_n)\le t_n\}\right)\le \frac{1}{C}S \left(y,\{z:\,S(z,x_n)\le t_n\}\right)\leq \frac{\d_n}C
\end{align*}
where the latter inequality is due to the subadditivity of $S$.
Therefore, by \eqref{eq7} and \eqref{eq8}, we deduce
\begin{align*}
  u(x_0,t_0)
  &\ge \min\{u_0(y):\,S(y,x_n)\le t_n+\d_n\}\\
&\ge \min\{u_0(y):\,\,d\left(y,\{z:\,S(z,x_n)\le t_n\}\right)\le \d_n/C\}\\
&\ge u (x_n,t_n)-\o(\d_n/C)
\end{align*}
where $\o$ is the modulus of continuity for $u_0$ in a   neighborhood of $x_0$.
This gives
\[u(x_0,t_0)\ge\limsup_{(x_n,t_n)\to (x_0,t_0)} u(x_n,t_n).\]
By $\{y \in \G:\,S(y,x_n)\le t_n\}\subset\{y \in \G:\,S(y,x_0)\le
t_0+\d_n\}$ we get in a similar way
\[u(x_0,t_0)\le\liminf_{(x_n,t_n)\to (x_0,t_0)} u(x_n,t_n).\]
We now prove that $u$ is a solution of \eqref{HJe}.
We only prove that $u$ is a supersolution at $(x_0,t_0)$ with $x_0=x_i \in V$, since the other cases can be proved  as in  the Euclidean case.
Assume by contradiction that there is an admissible test function $\phi$  such that $u-\phi$ has    a  local minimum  at $(x_0,t_0)$  with $\phi(x_0,t_0)=u(x_0,t_0)=\a$ and such that
\begin{equation}\label{contradiction}
 \max_{j\in Inc_i}\{\pd_t  \phi_j (x_0,t_0)+H_j (x_0, D_j\phi(x_0,t_0))\}\le -\d<0.
\end{equation}
Observe that
\begin{equation}\label{eq9}
\{(x,t):\,u(x,t)\le \a\}=\{(x,t):S(\{u_0\le \a\},x)\le t\}.
\end{equation}
Assume first that $S(\{u_0\le \a\},x_0)>0$ and  define $w(x,t)=S(\{u_0\le \a\},x)-t$. We claim that $\phi$ has a local maximum on the set $\{w\le 0\}$  at $(x_0,t_0)$. In fact if $(x,t)\in\{w\le 0\}$ then by \eqref{eq9}, $u(x,t)\le \a$  and since
\begin{equation}\label{eq10}
0=u(x_0,t_0)-\phi (x_0,t_0)\le u(x,t)-\phi (x,t)
\end{equation}
we get $\phi(x,t)\le  u(x,t)\le\a\le \phi(x_0,t_0)$ and the claim is proved.
By Proposition~\ref{prp:dist} and Lemma~\ref{caratstaz} $w$ is supersolution to \eqref{HJe} at $x_0$ and  therefore  Lemma~\ref{geometric}  gives a contradiction to \eqref{contradiction}.\par
If $S(\{u_0\le \a\},x_0)=0$, we claim that $(x_0,t_0)$ is a local maximum point for $u$. In fact, $S(\{u_0\le \a\},x_0)=0\le t_0-\eta$ for some $\eta>0$. If $(x,t)$ is such that $\max\{S(x,x_0),|t-t_0|\}\le \d/2$ with $\d<\eta$, then
\[S(\{u_0\le \a\},x)\le S(\{u_0\le \a\},x_0)+S(x_0,x)\le t_0-\eta+\d/2\le t\]
 hence $u(x,t)\le \a=u(x_0,t_0)$ and the claim is proved. By \eqref{eq10}
$(x_0,t_0)$ is also a local maximum point for $\phi $. By \eqref{H3} and \eqref{contradiction}, we get $\phi_t(x_0,t_0)<0$ and therefore a contradiction to $(x_0,t_0)$ being a local maximum point for $\phi $.
\end{Proofc}

\section{An application: the blocking problem}\label{sect:3}
In this section we provide a concrete application of our results: now, the network $\G$ represents  an oil pipeline (a network of computer,  the circulatory system, etc.) and at initial time a fire breaks up  in the region $R_0\subset \G$ (a virus is detected in a subnet, an embolus occurs in some vessel). The speed of propagation of  the fire is known but it may depend on the state variable (and, in particular, on the edge of the network). Our aim is to determine an  optimal strategy to stop the fire and to minimize the burnt region.\par
As in  the flame propagation model described in
\cite{Ba}, let $R_0$ be the
initial burnt region and $R_t$ the region burnt at time $t$. Assume that the front $\partial R_t$
propagates in the outward normal direction  to the front itself. Then  $R_t$ is given by the $0$-sublevel set of a viscosity solution of \eqref{HJe} -\eqref{IC} where the initial datum $u_0$ satisfies
  $R_0=\{x\in\G:u_0(x)\le 0\}$.\par
Recalling  the representation formula \eqref{31} we observe that
 the  0-sublevel set of the solution of \eqref{HJe}-\eqref{IC} is
 given by
\begin{equation*}
 R_t=\{x\in\G:\, S(R_0,x)\le t\}
\end{equation*}
where $S$ is defined as \eqref{dist}. Note that, since $\G$ is composed by a finite  number of bounded edges and therefore its total length is finite, then the burnt region
\[ R =\{x\in\G:\, S(R_0,x)<\infty\}=\cup_{t\ge 0} R_t\]
coincides with $\G$. In other words, without any external intervention, the pipeline will be completely burnt in a finite time.\par
We assume that an operator, located at $x_0\in V$ (the ``operation center''), can block the fire by closing the junctions of the pipeline (i.e., vertices of the network) and that this operation is effective only after a delay which depends on the distance of the junction from $x_0$.\par
Our problem is reminiscent of other models described in literature (for instance, see \cite{kb, vok} and references therein) which concern the control of some diffusion in a network (e.g. minimizing the spread of a virus or maximizing the spread of an information).
In this framework, let us stress the main novelties of our setting: in our model, the diffusion has positive finite speed and it affects both vertices and edges, the spread is not reversible (namely, ``infected'' points cannot become again ``healthy'') and the effect of the operator's action has finite speed (in other words, it is effective after a delay depending on the distance from the operation center).
\begin{Definition}
An admissible strategy $\s$ is a subset of $V$ such that
\begin{equation}\label{admissible}
S (R_0,x_i) \ge \delta d (x_0,x_i)\qquad\forall x_i\in \s
\end{equation}
where $\d$ is a given nonnegative  constant.
We denote by $V_{ad}$ the set of the vertices which satisfy the admissibility condition  \eqref{admissible} and by  $\gS_{ad}$ the set of the admissible strategies.
\end{Definition}
\begin{Remark}
Condition \eqref{admissible} means that the time to reach the  vertex $x_i\in \s$   from $x_0$ at the velocity $1/\d$ is less than or equal to the time the  fire front reaches $x_i$. Therefore the junction $x_i$ can be blocked   before the front goes  through it.
\end{Remark}
Given a  strategy  $\s\in \gS_{ad}$, we denote $S^\s:\G\times \G\to [0,\infty]$ the distance restricted to the trajectories not going
 through a vertex in   $\s$,  i.e.
\[
   S^\s(y,x):=\inf\left\{\int_0^t s (\xi(r),\dot\xi(r))dr:\,t>0,\,\xi\in B^t_{y,x}\,\text{s.t. $\xi(r)\not\in \s$ $\forall  r\in [0,t]$}\right\}
\]
with $S^\s(y,x)=\infty$ if there is no admissible curve joining $y$ to $x$. We also set
\begin{align*}
& R^\s_t:=\{x\in\G:\, S^\s(R_0,x)\le t\}\\
& R^\s :=\cup_{t\ge 0} R^\s_t=\{x\in\G:\, S^\s(R_0,x)<\infty\}
\end{align*}
which are respectively the region burnt at time $t$ and the total burnt region using the strategy $\s$.
Observe that
\begin{itemize}
  \item  if $\d$ is very small, then the optimal strategy is given  by the endpoints of the edges containing $R_0$;
\item if $\d$ is very large and $x_0\in R_0$, then every strategy is useless since the whole pipeline will burn whatever the operator does.
\end{itemize}
Aside  the previous  simple cases an optimal strategy for the blocking problem may be not obvious and we aim to find an efficient way to compute it.
To find a strategy which minimizes the burnt region, we first
 give  a characterization of   $R^\sigma$     in terms  of a problem satisfied by the distance $S^\s(R_0,\cdot)$.
\begin{Proposition}
Given $\s\in \gS_{ad}$, set $u(x)=S^\s(R_0,x)$ and $\cR=R^\s $. Then
 \begin{itemize}
 \item[i)] $u\in C^0(\cR)$ and $u= +\infty$ in $\G\setminus \cR$. Moreover if $x_i\in\s$ and  $j\in Inc_i$ is such that $e_j\subset  \cR$, then  $\lim_{x\to x_i,\, x\in e_j} u(x)=u_j(x_i)<\infty$;
 \item[ii)] $u$ is a viscosity solution of the  problem
\begin{equation}\label{HJvinc}
 \left\{
 \begin{array}{ll}
   H(x,Du)=1,\quad & x\in \cR\setminus(\ R_0\cup\s) \\
   u=0,  & x\in R_0;
 \end{array}
 \right.
\end{equation}
\item[iii)]  let $w\in C(\G)$ be  such that, defined $\cR_w=\{x\in\G:\, w(x)<\infty\}$,
\[
\left\{
\begin{array}{ll}
   H(x,Dw)\le1,\quad & x\in \cR_w\setminus  R_0 \\
   w=0,  & x\in R_0,
 \end{array} \right.
 \]
then $\cR\subset\cR_w$ and $w\le u$ in $\G$.
\end{itemize}
\end{Proposition}
\begin{Proof}
Note that, for $j\in J$,  either $e_j\subset  \cR$ or $  e_j\cap\cR=\emptyset$, i.e. an edge is either completely burnt or it cannot be reached by the fire.
The function $u$ can be discontinuous at $x_i\in\s$ and
\begin{itemize}
  \item if $x_i \in V\setminus \s$, then either $u_j(x_i)=\infty$ for all $j\in Inc_i$ if $x_i\in\G\setminus\cR$ or  $u_j(x_i)<\infty$ for all $j\in Inc_i$ if $x_i\in \cR$;
  \item  if $x_i \in \s$, then  either $u_j(x_i)=\infty$ for all $j\in Inc_i$  if $x_i\in\G\setminus\cR$  or there exists $j\in Inc_i$ such that $u_j(x_i)<\infty$ if $x_i\in \cR$   and in this case $u_j(x_i)=\sup_{e_j} u_j$.
\end{itemize}
Actually, if $x_i\in\s$ and   $u_j(x_i)<\infty$,  an admissible trajectory for $S^\s$ connecting $x_i$ to $R_0$ and containing the edge $e_j$, $j\in Inc_i$, necessarily  enters from $x_i$ into $e_j$. Hence  $u(x)$ is increasing for $x\in e_j$, $x\to x_i$ and $\lim_{x\in e_j,\,x\to x_i } u_j(x)=u_j(x_i)$. \par
In $\cR\setminus \s$,  $S^\s$ locally behaves as the distance $S$ defined  in \eqref{dist}. Therefore the continuity of $u$ in $\cR\setminus\s$ and the sub- and supersolution properties in the open set $\cR\setminus(\ R_0\cup\s)$ are obtained by the same arguments of \cite[Prop.4.1]{sc}.\par
To prove {\it iii)}, assume by contradiction that there exists  $x_0\in\cR$ such that $u(x_0)<w(x_0)$.
For any $x,y\in \cR$ such that $S^\s(x,y)<\infty$, a minimizing trajectory always exists  since (up to reparametrization) there is only  a finite number of trajectories connecting the two points.\par
Hence let $\xi$ be an admissible curve for $S^\s$ such that $\xi(0)=y_0\in R_0$, $\xi(T)=x_0$ and $u(x_0)=\int_0^Ts(\xi(r),\dot\xi(r))ds $.
Let $t_0=0<t_1<\dots<t_{M+1}=T$ such that, for any $m=0,\dots,M$,   $\xi([t_m,t_{m+1}])\subset  e_{j_m}$ for some $j_m\in J$, $\xi(t_{i_m})=x_{i_m}\in V$ and $\pi_{j_m}^{-1}\circ \xi\in C^1(t_m,t_{m+1})$. Clearly   $u(\xi (t))=\int_0^t s(\xi(r),\dot\xi(r))dr$ for $t\in [0,T]$.\par
 If $w$ is a subsolution to \eqref{HJvinc}, then by the coercivity of $H$,  $w$ is Lipschitz continuous in $\cR_w\setminus R_0$ and therefore $H(x, Dw)\le 1$ a.e. on $\cR_w\setminus R_0$. Moreover, by the definition of the support function $s$,  we have $H(x,p)\le 1$ if and only if $\sup_{q\in\R}\{pq-s(x,q)\}\le 0$. Hence
\[a_{i_mj_m}\dot\xi(r) Dw(\xi(r))\le s(\xi(r),\dot\xi(r))\qquad \text{$\forall   r\in [t_m,t_{m+1}]$, $m=0,\dots,M$.}
\]
(the term $a_{i_mj_m}$ takes into account the orientation of the arc $e_{j_m}$).
Integrating the previous relation in $[0,t_1]$ and recalling that $u(y_0)=w(y_0)=0$, we get
\[w(\xi(t_1))\le \int_0^{t_1} s(\xi(r),\dot\xi(r))dr=u(\xi(t_1)).\]
 Iterating the same argument  in $[t_{m},t_{m+1}]$ we finally get $w(\xi(T))\le u(\xi(T))$ and therefore a contradiction since $\xi(T)=x_0$.  We conclude that $x_0\in\cR_w$ and $w(x_0)\le u(x_0)$.
\end{Proof}
We  now show that the strategy composed by all  the admissible nodes which are adjacent to a non admissible node is optimal, in the sense that it maximizes the  preserved region.
\begin{Proposition}\label{optbl}
The admissible strategy
\begin{equation*}
    \s_{opt}=\{x_i\in V_{ad}:\, \exists x_j\in V\setminus V_{ad},\, e_k\in E\,\text{s.t.}\,x_i,\,x_j\in\overline e_k\}
\end{equation*}
satisfies: for any $\s\in \Sigma_{ad}$, $\G\setminus R^{\s}\subset \G\setminus R^{ \s_{opt}}$.
\end{Proposition}
\begin{Proof}
Assume by contradiction that there exist $\s\in  \Sigma_{ad}$ and  $x_0\in(\G\setminus R^\s)\cap R^{\s_{opt}}$. Hence there exists an admissible trajectory $\xi$ for $\s_{opt}$ connecting $x_0$ to $R_0$, i.e. there exists $y_0\in R_0$, $t>0$ and $\xi\in B^t_{y_0,x_0}$ such that  $\xi(r)\not\in \s_{opt}$ for $r\in [0,t]$. Note that $\s_{opt}$ disconnect the subgraph containing the admissible
vertices $V_{ad}$ by the one containing the non admissible vertices $V\setminus V_{ad}$ and therefore $\xi([0,t])$
 is contained in the subgraph with vertices $(V\setminus V_{ad})\cup \s_{opt}$. Since $\s\subset V_{ad}$,  then $\xi$ is also admissible for $S^{\s}(y_0,x_0)$ and therefore a contradiction to $x_0\in \G\setminus R^\s$.
\end{Proof}

\begin{Remark}
It is possible to consider a cost functional on the set of the admissible strategies $\G_{ad}$ which takes into account not only the part of the  network destroyed by the fire but also other terms, such as  the cost of blocking a given junction.
Consider  the cost functional $ { \cal I}:\G_{ad}\to\R$  given by
\[ {\cal I}(\s)=\sum_{x_i\in \s} \a_i  +\sum_{e_j\subset R^\s } \b_j\]
(recall that either $e_j\subset  R^\s $ or $  e_j\cap R^\s =\emptyset$).
The first term represent the cost of blocking the node $x_i$ and can depend on various parameter (distance of the
node from $x_0$, accessibility of $x_i$, cost of blocking $x_i$,  etc.) while the second term is the cost of the burnt region with a given cost $\b_j$ for each arc.
Clearly, the minimum of   ${ \cal I}(\cdot)$  exists since $\gS_{ad}$ is finite, but it seems more difficult to characterize the optimal strategy.
\end{Remark}
%
 \subsection{Numerical simulations}\label{S5}
In this section we propose a numerical method to compute the optimal strategy for the blocking problem. The scheme is based on a finite difference approximation of the stationary problem \eqref{HJvinc}; for simplicity,  we only consider the case of the eikonal Hamiltonian $H(x,p)=|p|/{c(x)}$.\par
On each interval  $[0,l_j]$  parametrizing the arc $e_j$,  we consider  an uniform  partition $y_{j,m}=m  h_j$ with   $M_j= {l_j}/{h_j}\in \N$ and $m =0,\dots,M_j$. In this way   we obtain a   grid $\mathcal{G}^h=\{x_{j,m}=\pi_j(y_{j,m}),\; j\in J, \; m=0,\dots,M_j\}$   on the network $\Gamma$.
 We    define  $\mathcal{R}_0^h=\mathcal{G}^h\cap R_0$, the set of the nodes in the initial front.
For $x_1,x_2\in \mathcal{G}^h$, we say that $x_1$ and $x_2$ are adjacent and we write  $x_1\sim  x_2$ if and only if they are the image of two adjacent grid points, i.e. $x_k=\pi_j(y_k)$,  for $y_k\in [0,l_j]$,  $k=1,2$,   $j\in J$ and $|y_1-y_2|=h_j$. Note that if  $x_i\in V$ is a vertex of $\G$, then   the nodes of the grid $\mathcal{G}^h$ adjacent to $x_i$  may  belong to different arcs.\\
 We compute the optimal strategy by means of  the following Algorithm, based on the results of Prop. \ref{optbl}.\par
{\bf{Blocking strategy [B]}}
\begin{enumerate}
\item
In the first step we solve the  front propagation problem on the network computing
the  approximated time $u^h(x)$ at which a node $x_{j,m} \in \mathcal{G}^h$  gets burnt
 \begin{equation}{\label{HJscheme1}}
\left\{
\begin{array}{ll}
 \max_{x\in \mathcal{G}^h,\,  x \sim x_{j,m} }\left\{-\frac{1}{h_j} (u^h(x)-u^h(x_{j,m})) \right\}-c(x_{j,m})=0& {x_{j,m}}\in \mathcal{G}^h\\[6pt]
 u^h(x_{j,m})=0& x_{j,m}\in \mathcal{R}_0^h\\[6pt]
\end{array}
\right.
\end{equation}
Note that if $x_{j,m}$ coincides with a   vertex $x_i\in V$, the approximating equation reads as
\[\max_{j\in Inc_i}\,\max_{x \in \mathcal{G}^h\cap e_j,\, x \sim x_{j,m} }\left\{-\frac{1}{h_j} (u^h(x)-u^h(x_{j,m})) \right\}-c(x_{j,m})=0.\]
The discrete function  $u^h:\mathcal{G}^h \to\R$  is such that $u^h(x_{j,m})\simeq u(x_{j,m})$, where $u(x)= S(R_0,x)$.

\item
 In the second step we  determine the   vertices   which satisfy the admissibility condition \eqref{admissible}.
 We define $V^h_{ad}=\{x_i\in V : w^h(x_i)<u^h(x_i)\}$, where $w^h:\mathcal{G}^h \to\R$  represents the approximated time to reach a node $x\in\mathcal{G}^h$, starting from the operation center $x_0$   and moving with a constant  speed $1/\delta$. The   function $w^h$ is computed by means of the   finite difference scheme
\begin{equation}{\label{HJscheme2}}
\left\{
\begin{array}{ll}
 \max_{x\in \mathcal{G}^h,\,x \sim x_{j,m} }\left\{-\frac{1}{h_j} (w^h(x)-w^h(x_{j,m})) \right\}-\frac{1}{\delta}=0& x_{j,m}\in \mathcal{G}^h\setminus\{x_0\}\\[6pt]
 w^h(x_0)=0.&   \\[6pt]
\end{array}
\right.
\end{equation}
\item We define  the approximated optimal strategy  by setting
$$\s^h_{opt}=\{x_i\in V^h_{ad}:\, \exists x_j\in V\setminus V^h_{ad},\, e_k\in E\,\text{s.t.}\,x_i,\,x_j\in e_k\}
$$
and, for $\s=\s^h_{opt}$,  we   compute the corresponding  approximate distance    by solving the following finite difference scheme
\begin{equation*}
\left\{
\begin{array}{ll}
 \underset{{x\in  \mathcal{G}^h,\,x \sim x_{j,m} }}\max\left\{-\frac{1}{h_j} (u^h_\s(x)-u^h_\s(x_{j,m})) \right\}-c(x_{j,m})=0& x_{j,m}\in \mathcal{G}^h \setminus(\mathcal{R}_0^h\cup\s)\\[6pt]
 u^h_\s(x_{j,m})=0& x_{j,m}\in \mathcal{R}_0^h\\[6pt]
 -\frac{1}{h_j} (u^h_{j,\s}(x)-u^h_{j,\s} (x_{j,m}))-c(x_{j,m})=0 &\text{$x_{j,m} \in \sigma$,   $x\in \mathcal{G}^h$, $x\sim x_{j,m}$}
\end{array}
\right.
\end{equation*}
The discrete function  $u^h:\mathcal{G}^h \to\R$  is such that $u^h_\s( x_{j,m}) \simeq u(x_{j,m})$, where $u(x)=S^\s(R_0,x)$ solves $\eqref{HJvinc}$.
Note that as in the continuous case, the value of $u^h$ at $x_i\in \s$ can depend on the edge $e_j$ and in general the function is discontinuous at these points.
\end{enumerate}
\begin{Remark}
In this paper we do not analyze the properties  of the previous finite difference schemes. In any case, at least for    \eqref{HJscheme1} and \eqref{HJscheme2}, the well-posedness and the convergence of the schemes can be studied by adapting the techniques  in \cite{cfs}.
\end{Remark}
\subsection{Example 1: a simple network}
We consider a  network with a simple structure where the fire starts in   one vertex, $R_0=\{(0,0)\}$  and propagates with   speed $c=1$. \par
\begin{figure}[h!]
\begin{center}
\epsfig{figure=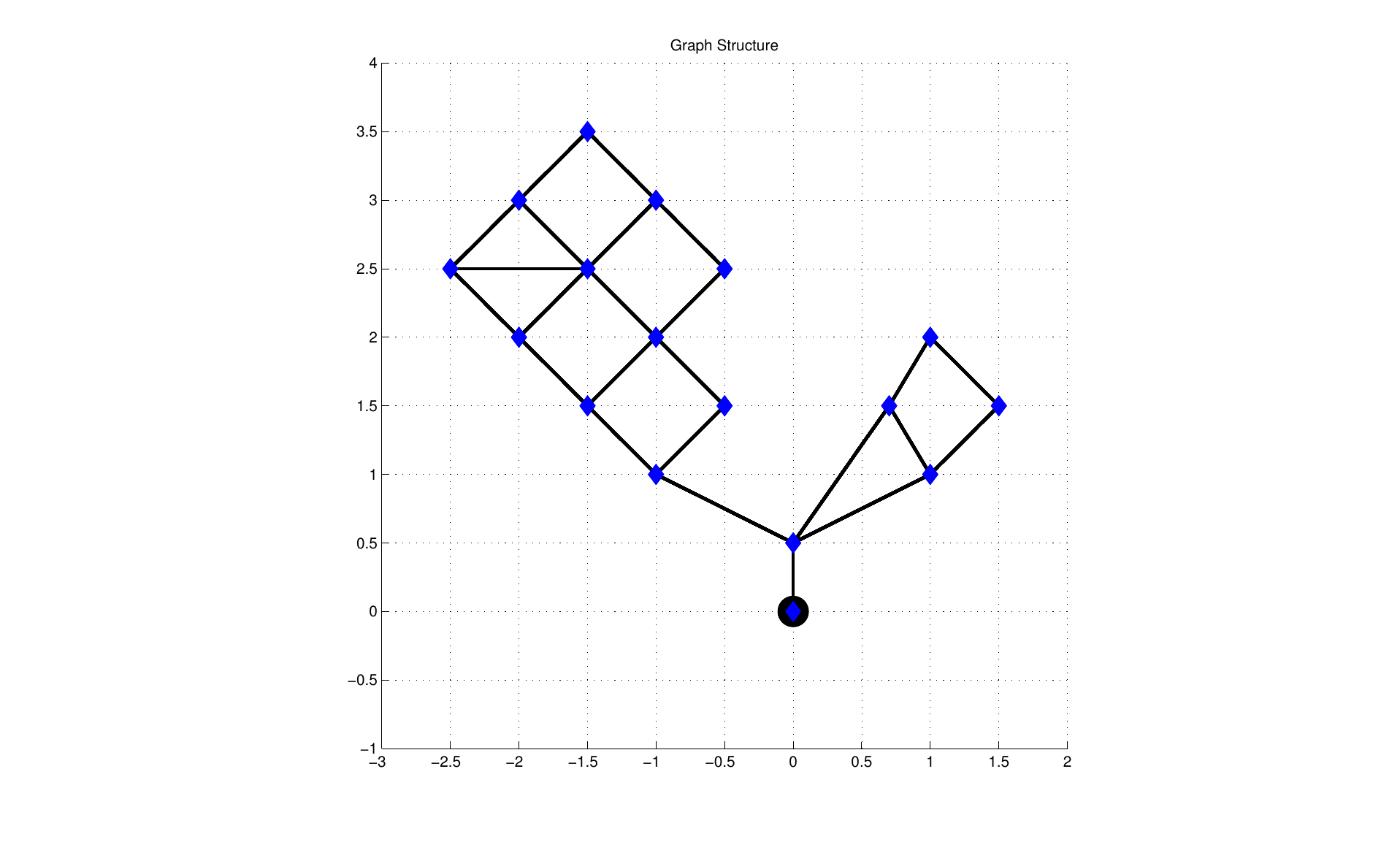,width=6.5cm}\hspace{-0.5cm}\epsfig{figure=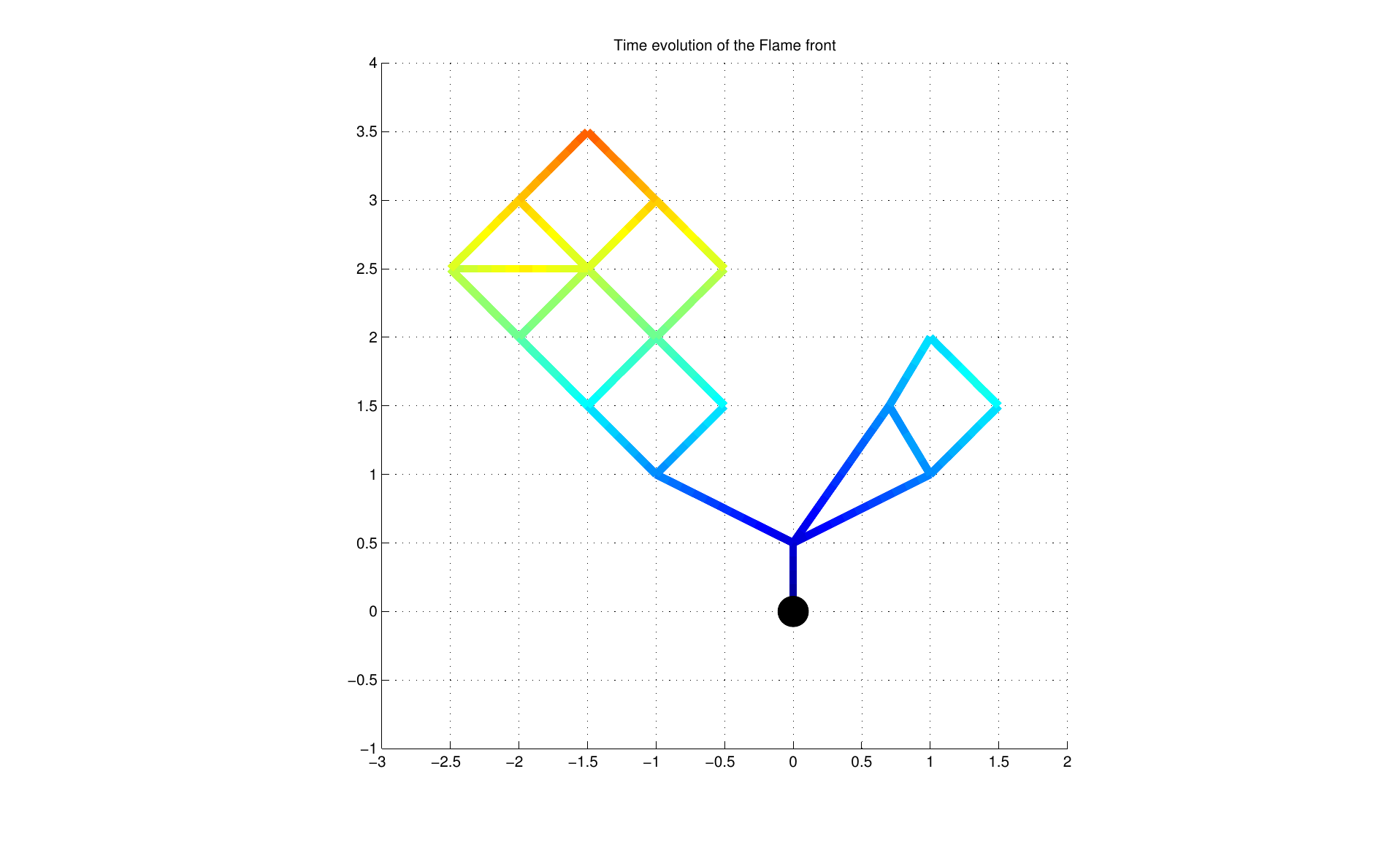,width=6.5cm}\hspace{-0.5cm}\epsfig{figure=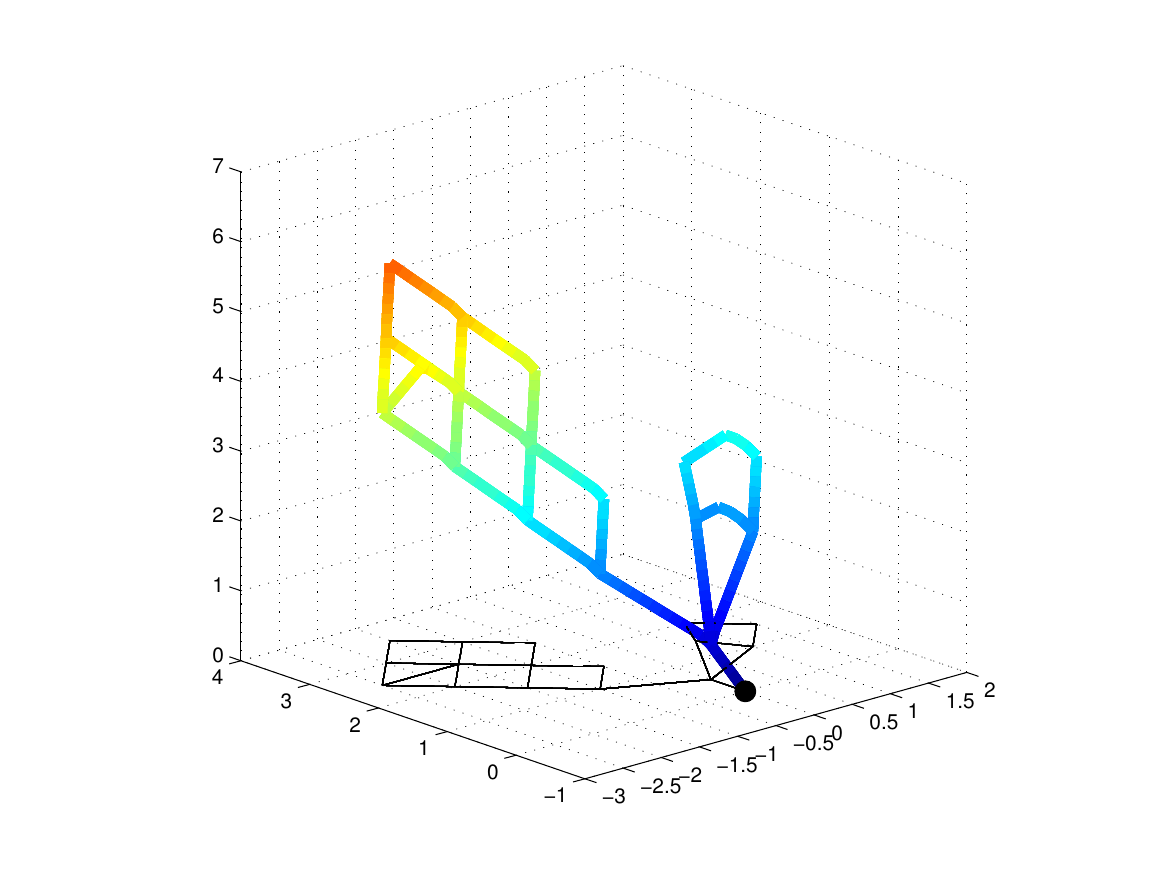,width=6.5cm}
\caption{Test1. Graph structure where $R_0$ is represented by the circle marker and the vertices by the rhombus markers (Top Left). Color map of  the  time $u_h(x)$ at which a node $x$ get burnt, computed by \eqref{HJscheme1}, (Top Right) and its 3D view (Bottom).}
\label{Test1grafo}
\end{center}
\end{figure}
\begin{figure}[h!]
\begin{center}
\epsfig{figure=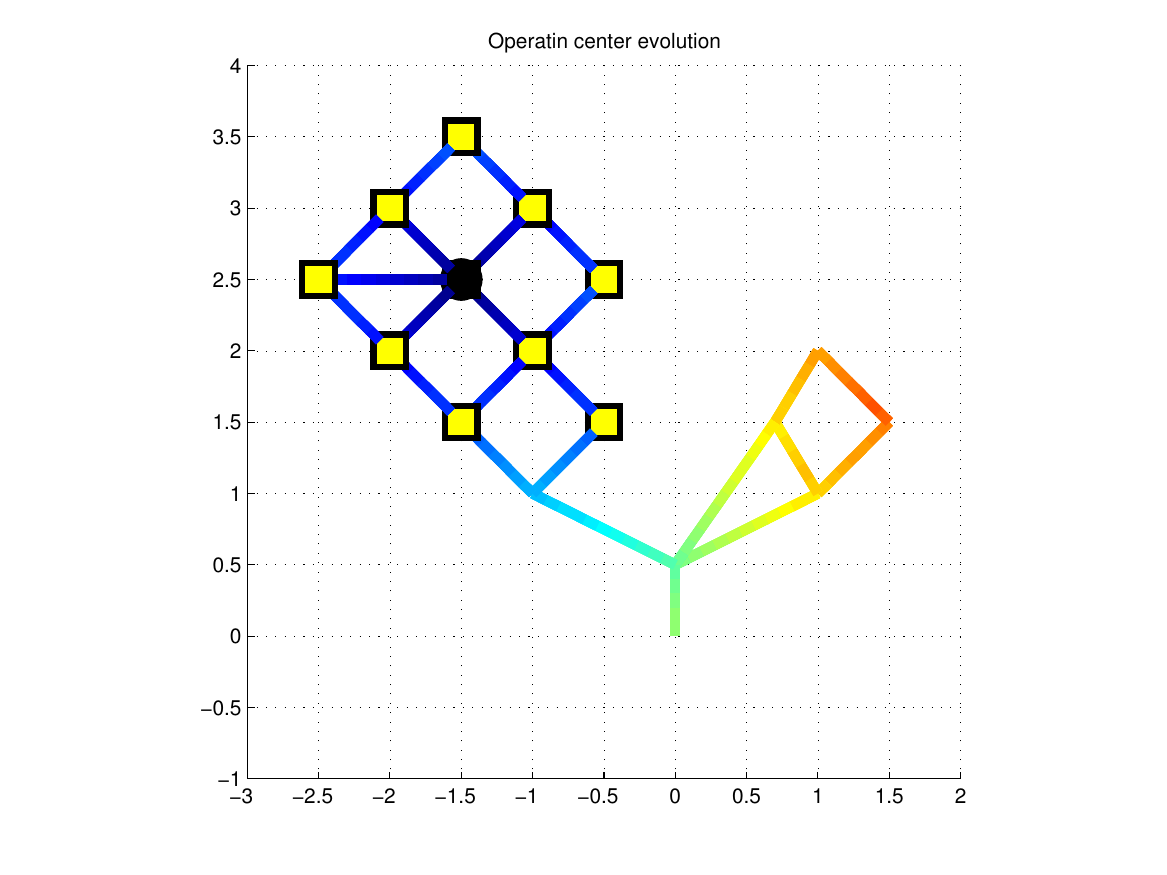,width=6.5cm}\epsfig{figure=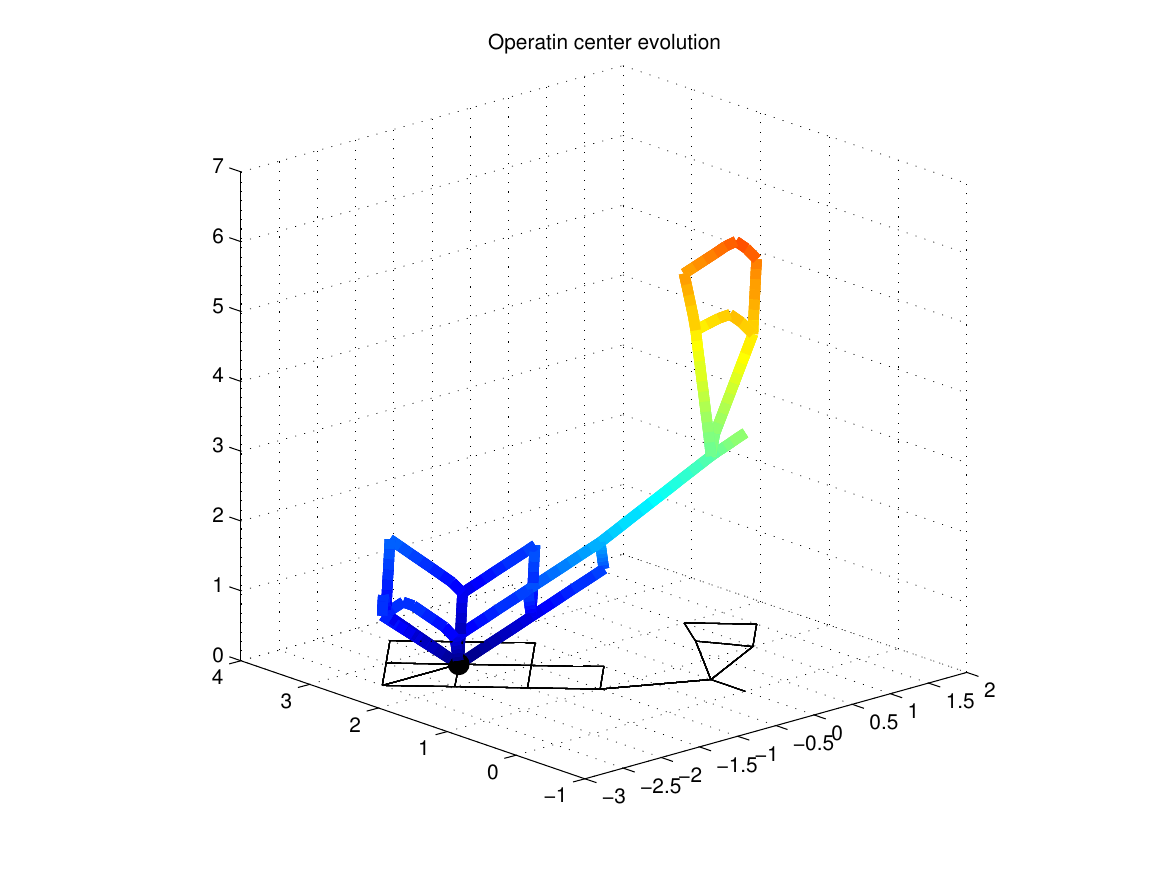,width=6.5cm}
\caption{Test1. Time to reach a point $x$ from the operation center $x_0$ (circle marker) and set of the admissible nodes $V^h_{ad}$ (square marker).
2D view (Left) and 3D view (Right).}
\label{Test1Vad}
\end{center}
\end{figure}

We first perform step $i)$  of Algorithm {\bf [B]} and we compute  the  approximated time $u^h(x)$ at which a node $x \in \mathcal{G}^h $  get burnt. The results are shown in Fig.\ref{Test1grafo} together with the graph structure.\\
Next, we perform step $ii)$ of Algorithm  {\bf [B]} . We suppose the operation center $x_0$ is located on the vertex $(-1.5,2.5)$  and the velocity to reach a node $x_i$ from $x_0$ is $\frac{1}{\delta}=1$. Using \eqref{HJscheme2}, we  compute the set of nodes $V^{h}_{ad}$. The result is shown in in Fig.\ref{Test1Vad}, the set of nodes  in $V^{h}_{ad}$ are represented by the square markers.\par
\begin{figure}[h!]
\begin{center}
\epsfig{figure=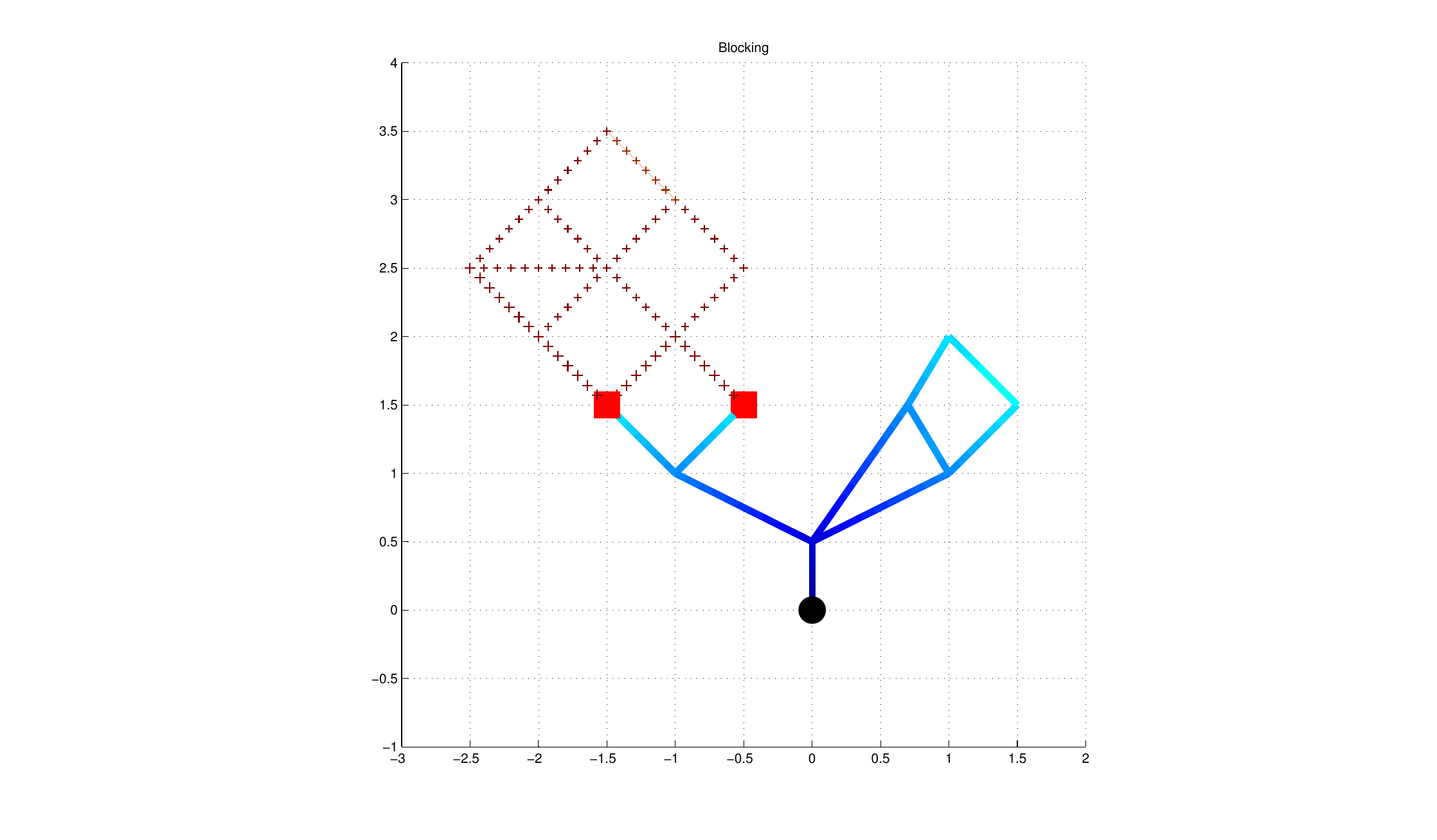,width=8cm}
\caption{Test1. Optimal blocking strategy $\s^h_{opt}$ (square marker), preserved network region (cross marker)  and minimum burnt network region (continuum line)  starting from $R_0$ (circle marker).}
\label{Test1blocking}
\end{center}
\end{figure}

Once computed  the set of admissible nodes $V^{h}_{ad}$, we can compute the  optimal strategy, following step $iii)$. The result is shown in Fig.\ref{Test1blocking}. It is clear, from the simple structure of the network,   that any other choice of $\s^h_{opt}$ would lead to a greater burnt region and consequently to a smaller preserved network region.\par
\subsection{Example 2: a more complex network}
\begin{figure}[h!]
\begin{center}
\epsfig{figure=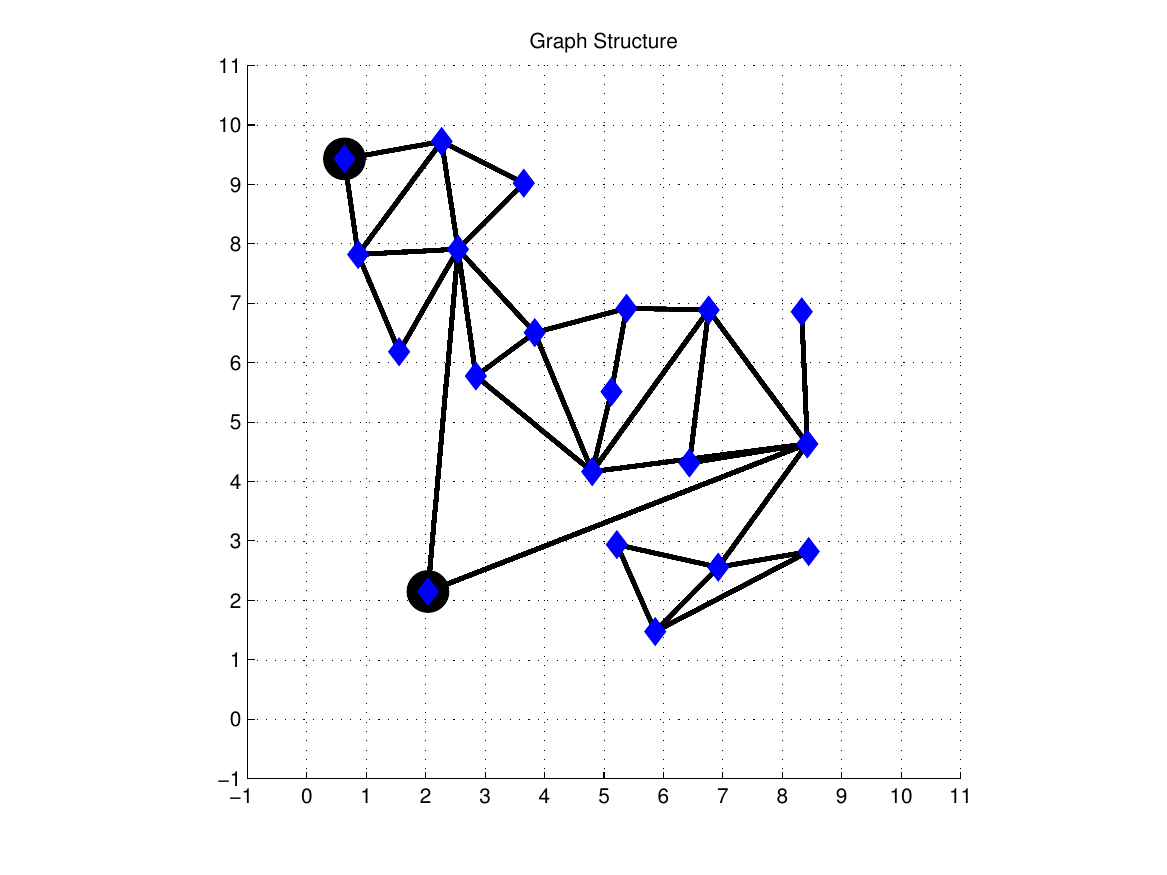,width=6.5cm}\hspace{-0.5cm}\epsfig{figure=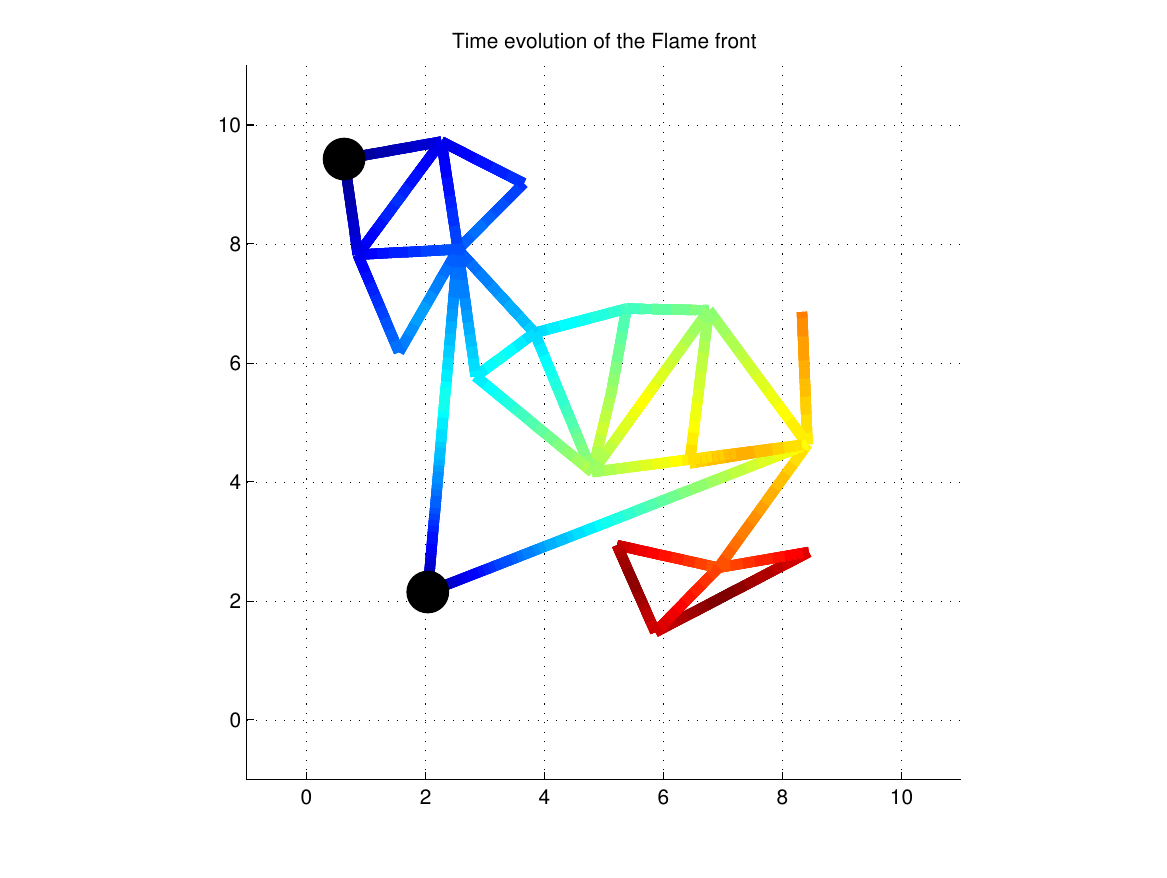,width=6.5cm}\hspace{-0.5cm}\epsfig{figure=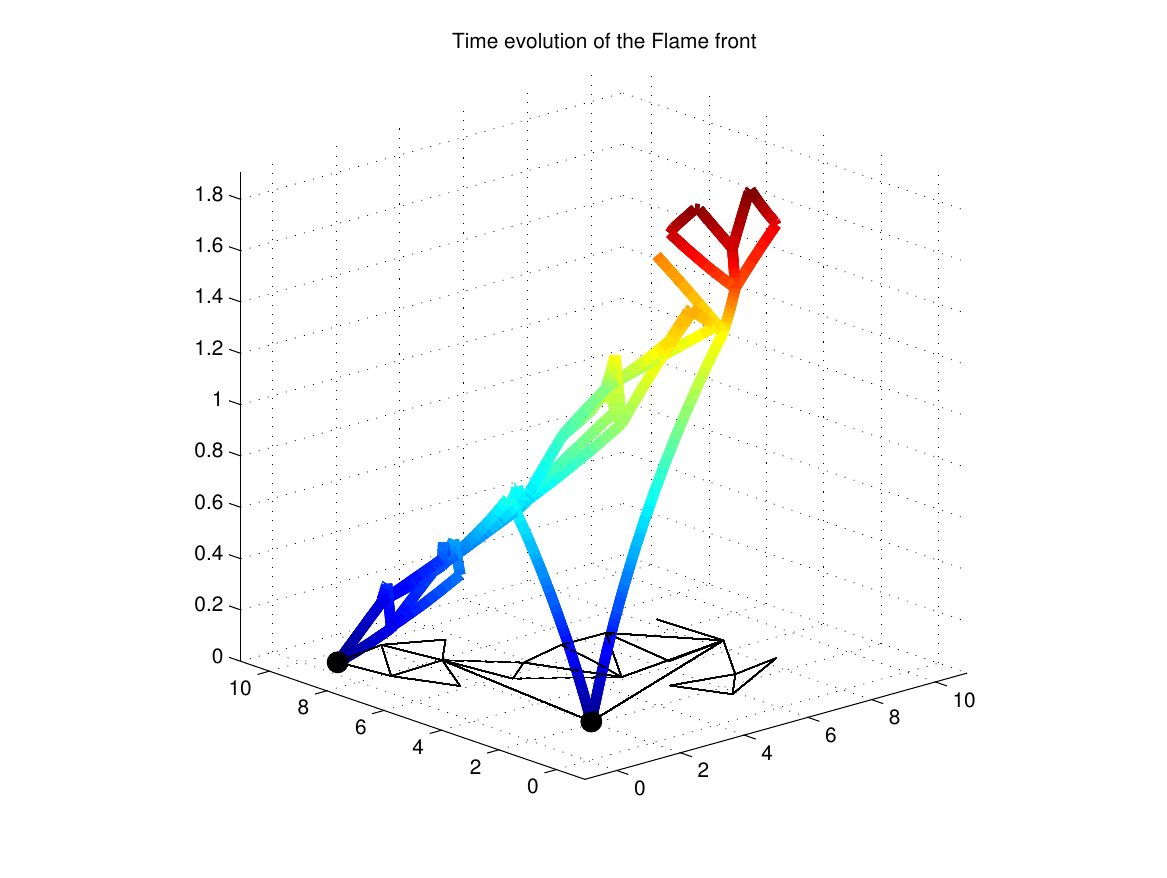,width=6.5cm}
\caption{Test2. Graph structure where $R_0$ is represented by the circle markers and the vertices by the rhombus markers (Top Left). Color map of  the  time $u_h(x)$ at which a node $x$ get burnt, computed by \eqref{HJscheme1} (Top Right), and its 3D view (Bottom).}
\label{Test2grafo}
\end{center}
\end{figure}

\begin{figure}[h!]
\begin{center}
\epsfig{figure=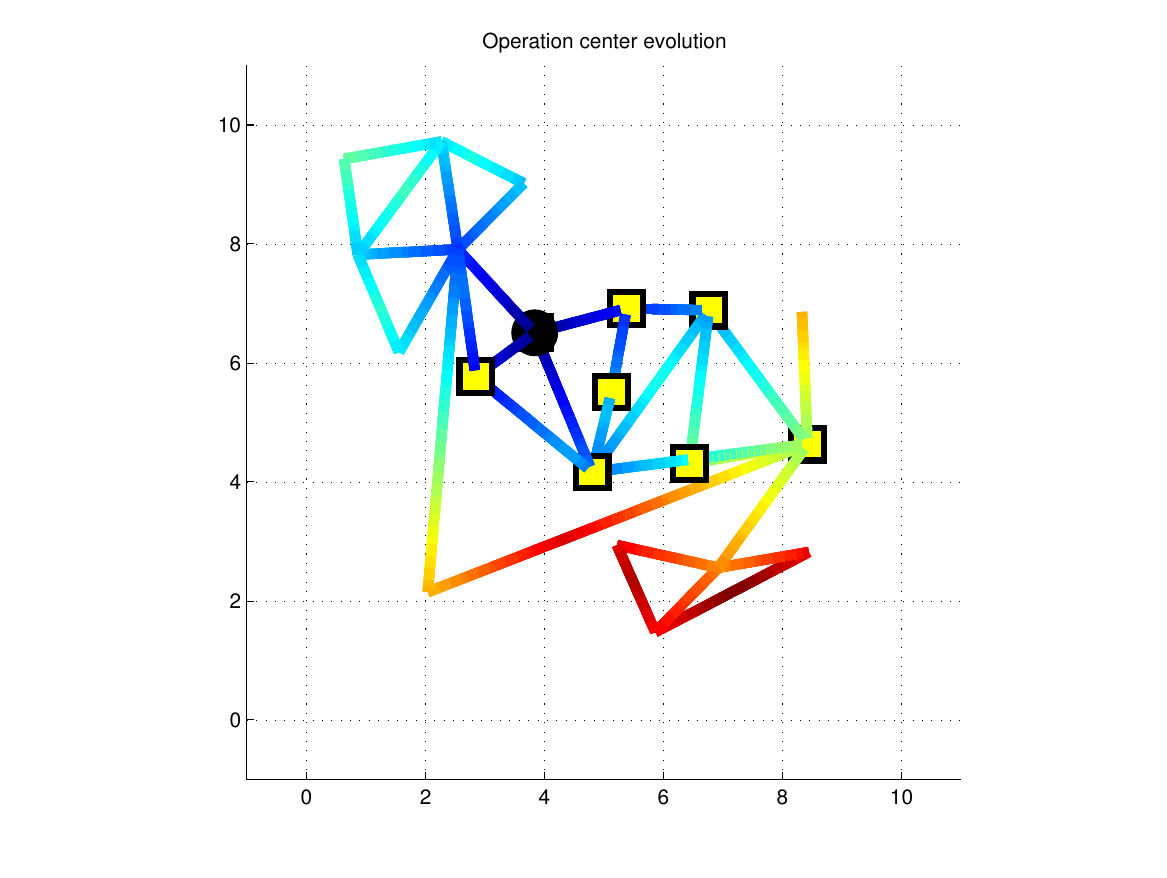,width=6.5cm}\epsfig{figure=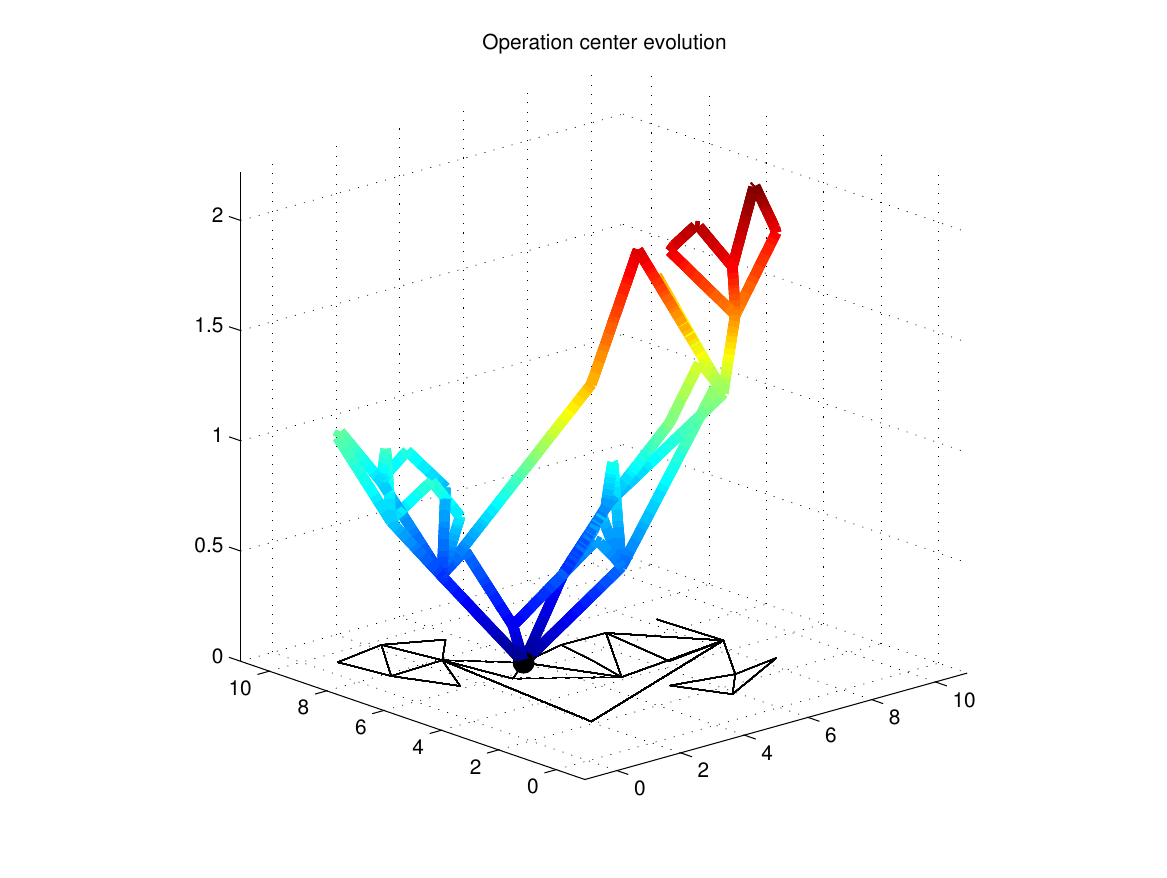,width=6.5cm}
\caption{Test2. Time to reach a point $x$ from the operation center $x_0$ (circle marker) and set of the admissible nodes $V^h_{ad}$ (square markers).
2D view(Left) and 3D view (Right).}
\label{Test2Vad}
\end{center}
\end{figure}

\begin{figure}[h!]
\begin{center}
\epsfig{figure=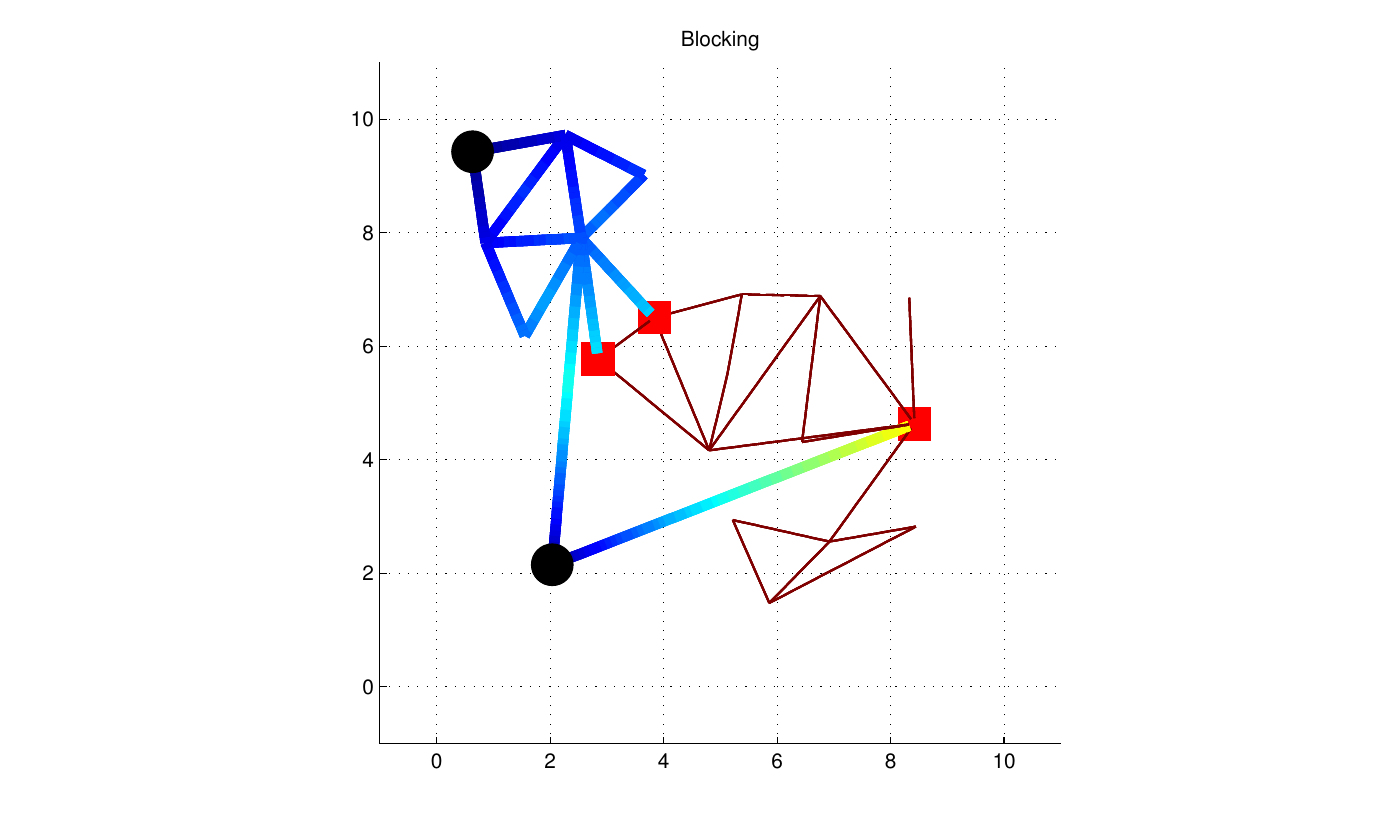,width=8cm}
\caption{Test2. Optimal blocking strategy $\s^h_{opt}$ (square markers), preserved network region (thin line)  and minimum burnt network region (thick line)  starting from $R_0$ (circle markers).}
\label{Test2blocking}
\end{center}
\end{figure}

\begin{figure}[h!]
\begin{center}
\epsfig{figure=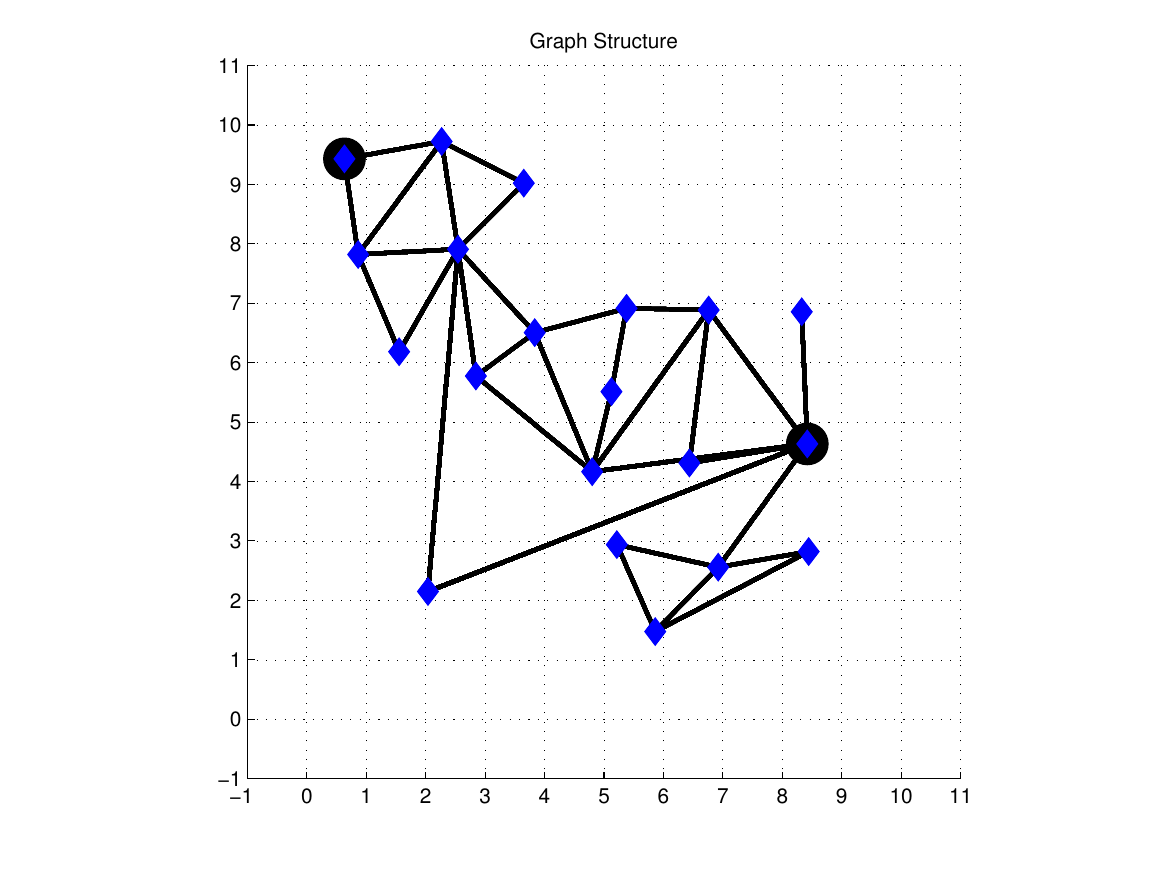,width=6.5cm}\hspace{-0.5cm}\epsfig{figure=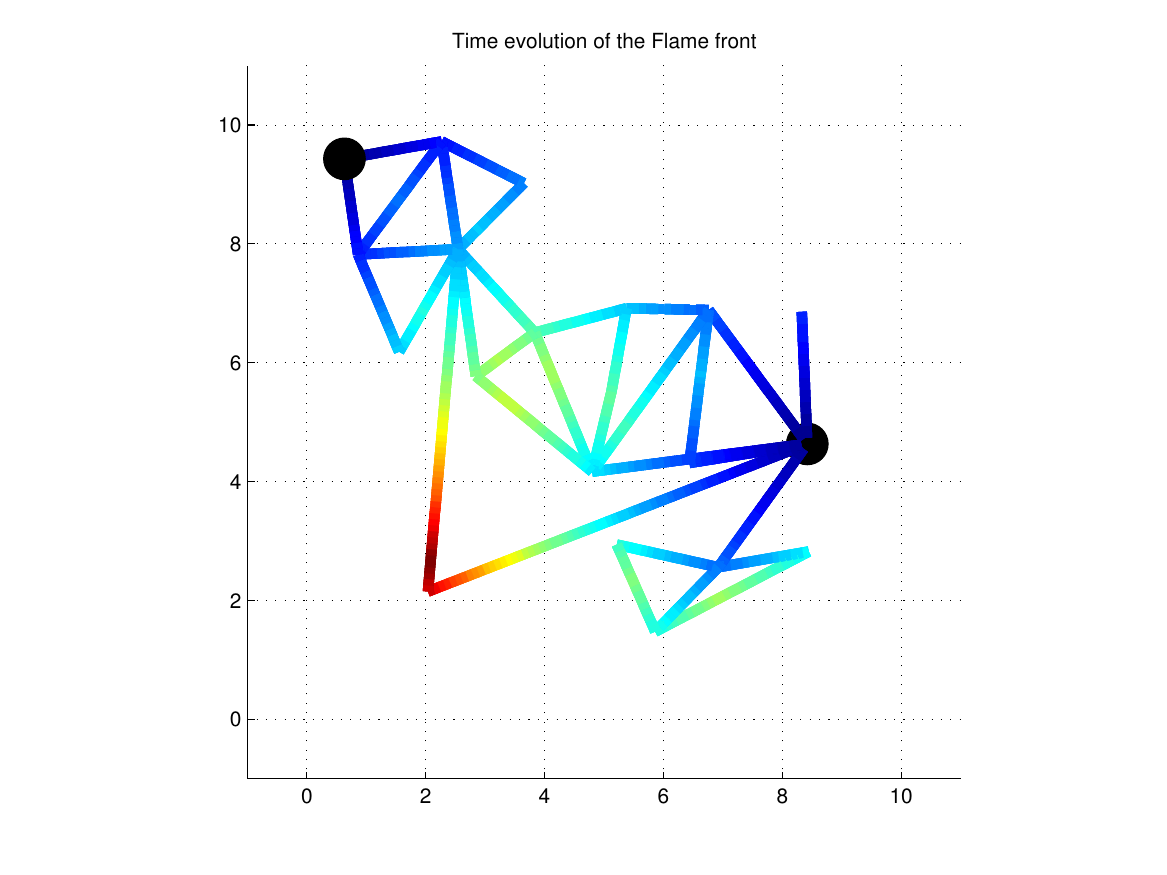,width=6.5cm}\hspace{-0.5cm}\epsfig{figure=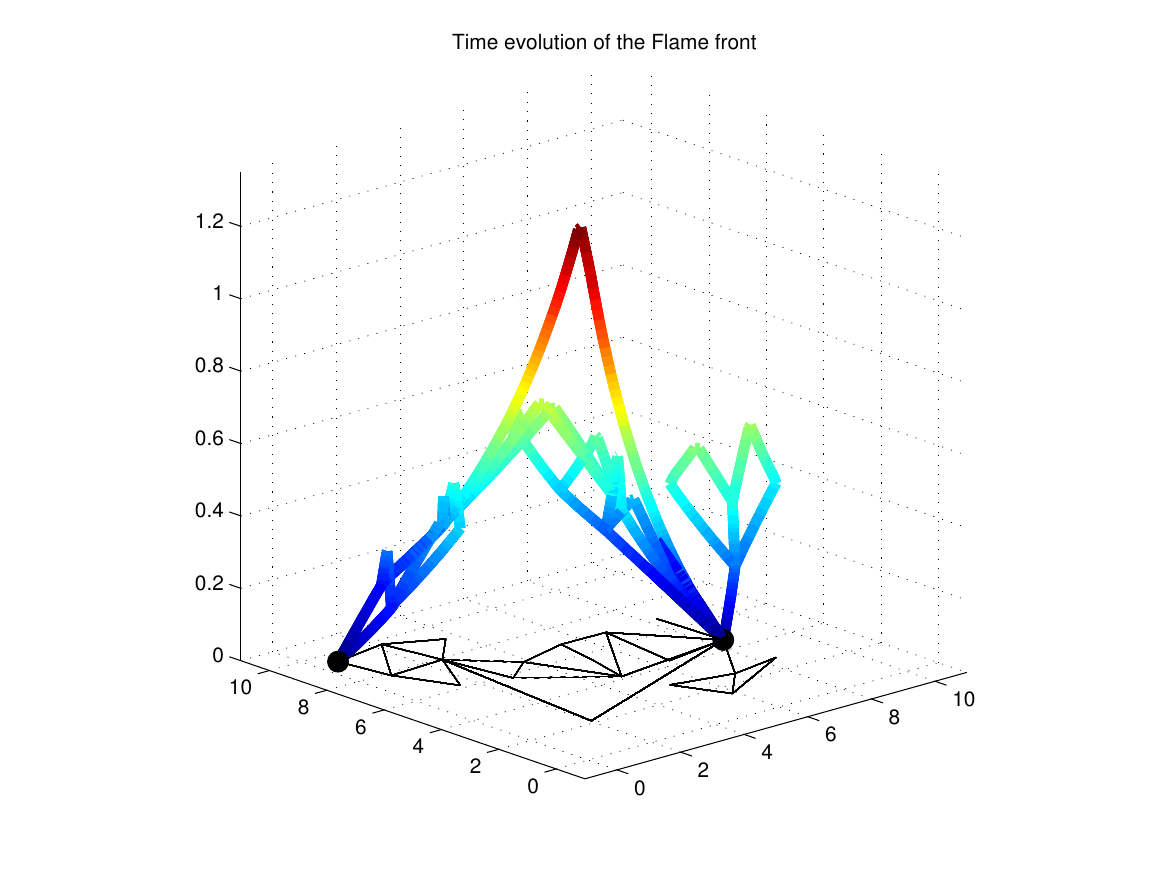,width=6.5cm}
\caption{Test3. Graph structure where $R_0$ is represented by the circle markers and the vertices by the rhombus markers (Top Left). Color map of  the  time $u_h(x)$ at which a node $x$ get burnt, computed by \eqref{HJscheme1}, (Top Right) and its 3D view (Bottom).}
\label{Test3grafo}
\end{center}
\end{figure}

\begin{figure}[h!]
\begin{center}
\epsfig{figure=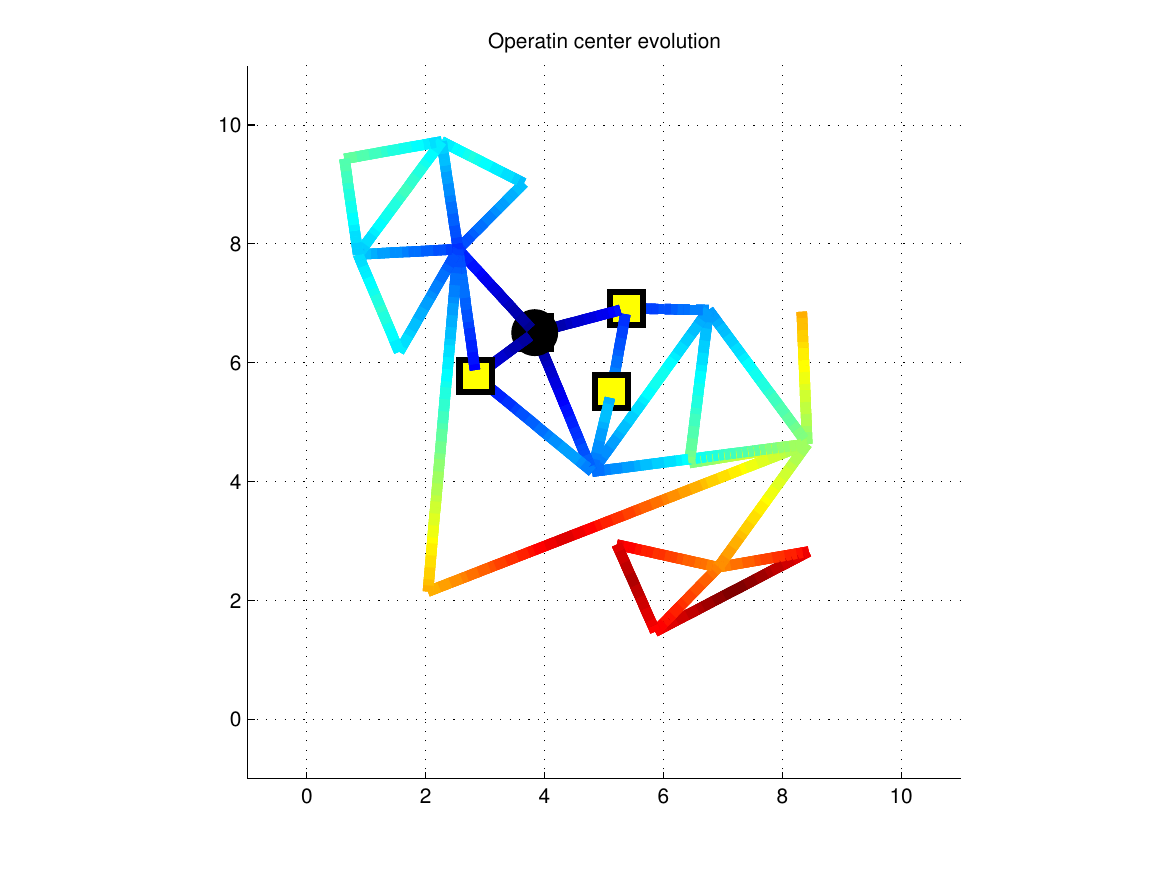,width=6.5cm}\epsfig{figure=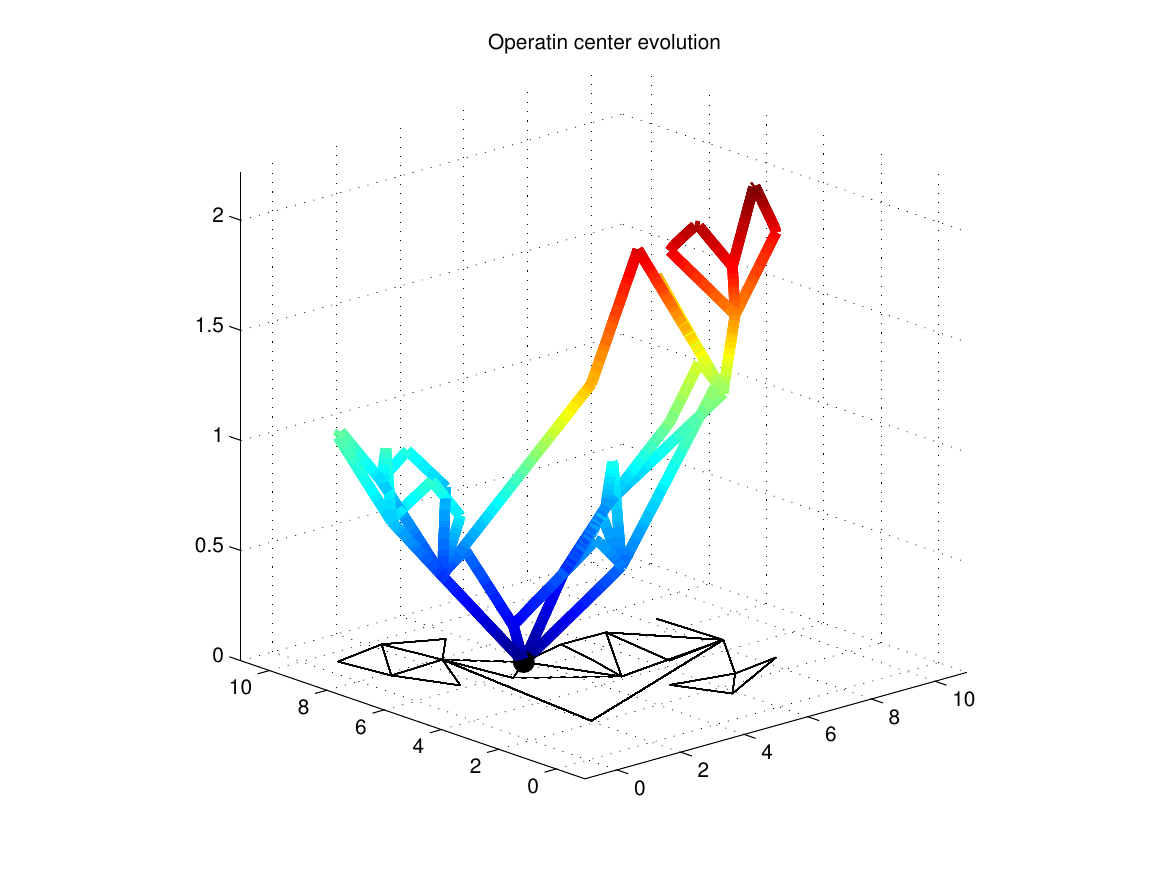,width=6.5cm}
\caption{Test3. Time to reach a point $x$ from the operation center $x_0$ (circle marker) and set of the admissible nodes $V^h_{ad}$ (square markers).
2D view(Left) and 3D view (Right).}
\label{Test3Vad}
\end{center}
\end{figure}

\begin{figure}[h!]
\begin{center}
\epsfig{figure=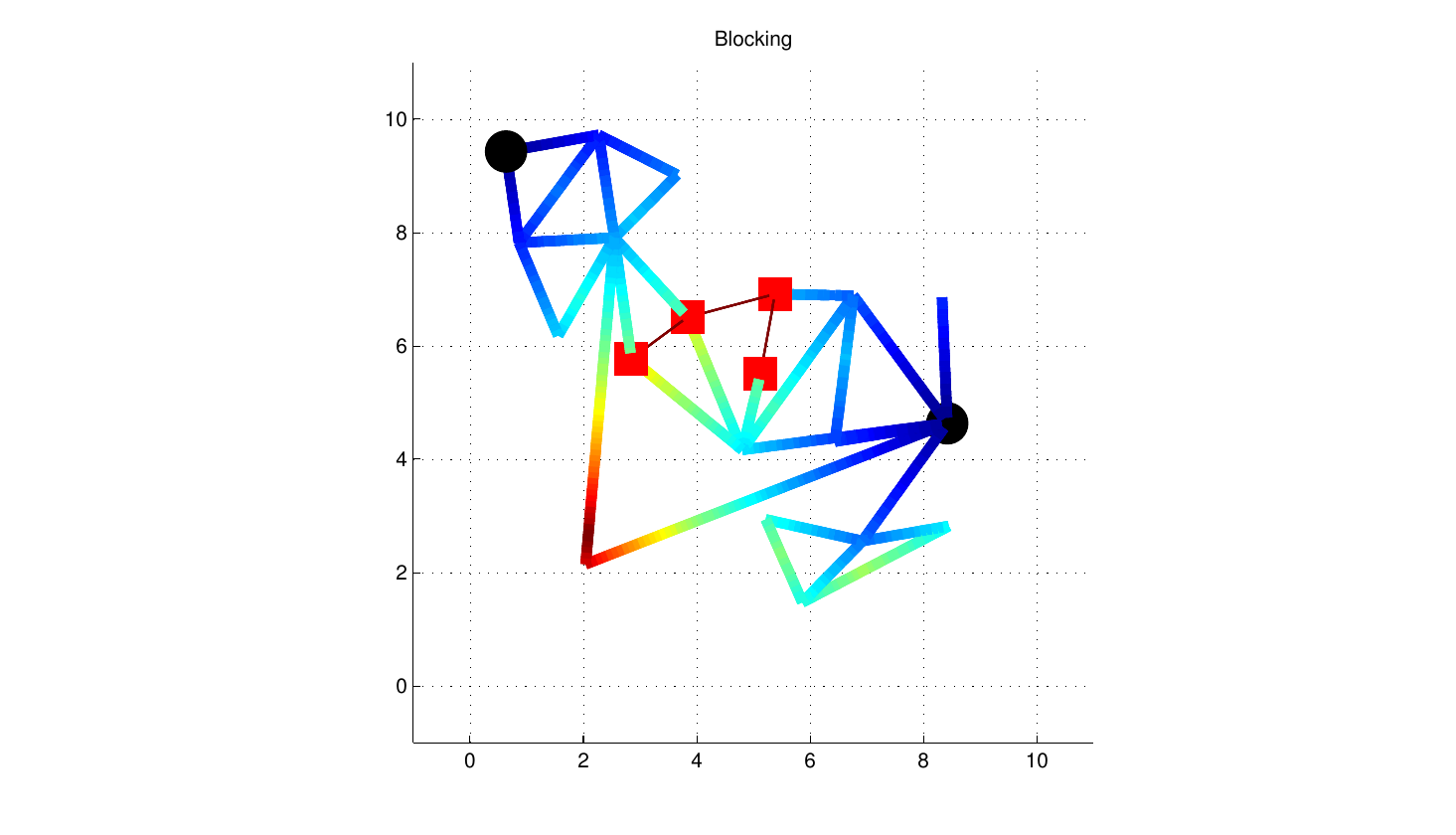,width=8cm}
\caption{Test3. Optimal blocking strategy $\s^h_{opt}$ (square marker), preserved network region (thin line)  and minimum burnt network region (thick line)  starting from $R_0$ (circle marker).}
\label{Test3blocking}
\end{center}
\end{figure}

We consider a more complex network, with  20 vertices and 32 arcs. We suppose the fire starts in two vertices and propagates with a non constant normal speed $c(x)=|x|$. \\
We proceed as in the first example and we compute the approximated time $u^h(x)$ at which a node $x \in \mathcal{G}^h $  get burnt. The results are shown in Fig.\ref{Test2grafo} together with the graph structure.\\
We suppose the operation center $x_0$ is located on the vertex $x_0=(3.8,6.5)$  and the velocity to reach a node $x_i$ from $x_0$ is $\frac{1}{\delta}=1/5$. Using \eqref{HJscheme2}, we  compute the set of nodes $V^{h}_{ad}$. The result is shown in in Fig. \ref{Test2Vad}, the set of nodes  in $V^{h}_{ad}$ are represented by the square markers.
Once computed  the set of admissible nodes $V^{h}_{ad}$, we can compute the  optimal strategy, following step $iii)$. The result is shown in Fig.\ref{Test2blocking}.
In this case we get the optimal solution  blocking only three vertices. By changing the set $R_0$ as in Fig.\ref{Test3grafo}, the region of the admissible node $V_{ad}^h$, shown  in Fig. \ref{Test3Vad}, turns out to be much smaller. In this case the optimal strategy is formed by all the  vertices in $V_{ad}^h$, and the preserved region becomes smaller than the previous case, see Fig. \ref{Test3blocking}.

\medskip
\noindent{\bf Acknowledgments} The third author is a member of Indam-Gnampa and he is partially supported by the Fond.-CaRiPaRo Project ``Nonlinear Partial Differential Equations: Asymptotic Problems and Mean-Field Games''.


\end{document}